\documentclass[10pt]{article}
\usepackage[utf8]{inputenc}
\usepackage{latexsym,amssymb,amsmath,url,authblk}

\def\N{\mathbb{N}}
\def\Q{\mathbb{Q}}
\def\Z{\mathbb{Z}}
\def\R{\mathbb{R}}
\def\C{\mathbb{C}}

\def\proof{\par\noindent{\em Proof. }}
\def\eproof{\hfill{$\Box$}\bigskip}

\def\ds{\dots}
\def\sus{\subset}
\def\al{\alpha}
\def\be{\beta}

\def\de{\delta}

\def\cc{\colon}
\def\ep{\varepsilon}

\newtheorem{thm}{Theorem}[section]
\newtheorem{prop}[thm]{Proposition}
\newtheorem{cor}[thm]{Corollary}

\newtheorem{lem}[thm]{Lemma}
\newtheorem{prob}[thm]{Problem}
\newtheorem{defi}[thm]{Definition}

\title{Transience of walks in the countable complete digraph with complex edge weights}
\author{M. Klazar\footnote{{\tt klazar@kam.mff.cuni.cz}}\ \ (UK, Praha) and R. Horsk\'y (V\v SE, Praha)}
\date{\today}

\begin{document}
\maketitle
\begin{abstract}
We extend a fragment of the symbolic method in enumerative 
combinatorics by allowing coefficients in generating functions to be defined by infinite series and generalize 
P\'olya's theorem on random walks. P\'olya's theorem determines limits of the probabilities that walks in the 
grid graph on $\Z^d$, starting at 
the origin $\overline{0}$, visit the given vertex $\overline{v}$. 
We generalize the grid graph to the countable complete directed graph 
$\vec{K}_{\N}$ on the vertices $\N=\{1,2,\ds\}$ and weight its edges by 
complex numbers. We obtain 16 generalizations of P\'olya's 
theorem. Our approach, in fact, generalizes Markov chains on $\N$ starting at $1$.
\end{abstract}

\tableofcontents

\section{Introduction}\label{sec_intro}

In this section, we recall the generating functions proof of P\'olya's theorem for 
$\overline{v}=\overline{0}$, briefly treat the case 
$\overline{v}\ne\overline{0}$, explain our generalization, and 
conclude with the plan of the 
article and a~discussion of follow-up works.\footnote{The previous 
version of this article, available as version~9 of \cite{arXivVers}, contained an introduction to the 
symbolic method in enumerative combinatorics. The readers of {\em 
S\'eminaire Lotharingien de Combinatoire} do not need it.}

\subsection{P\'olya's theorem for $\overline{v}=\overline{0}$}\label{subsec_PTvSMv0}

Let $d,n\in\N$ ($=\{1,2,\ds\}$) and $\overline{v}\in\Z^d$, where $\Z=
\{\ds,-1,0,1,\ds\}$ are the integers. Let $P_d(\overline{v},n)$ ($\in[0,1]$) 
be the probability that a~walk of length $n$ in the grid
graph on $\Z^d$, starting at the origin 
$\overline{0}=\langle0,0,\ds,0\rangle$ ($\in\Z^d$), later visits 
the vertex~$\overline{v}$. In 1921, G.~P\'olya proved a~theorem \cite{poly}, which we generalize in this article in several ways, according to which 
\begin{equation}\label{polya}
\lim_{n\to\infty}P_d(\overline{v},\,n)\left\{
\begin{array}{lll}
 =1 & \ds  & \text{if $d\le2$ and} \\
 \in(0,\,1) &  \ds & \text{if $d\ge3$}\,.
\end{array}
\right.\end{equation}
P\'olya's theorem is usually stated only for $\overline{v}=\overline{0}$, but he proved it for every 
vertex~$\overline{v}$, see the original article \cite{poly}. 
We noticed a~structural difference between the cases $\overline{v}=\overline{0}$ and 
$\overline{v}\ne\overline{0}$. For $\overline{v}=\overline{0}$, both cases $d\le2$ and $d\ge3$ in 
the dichotomy (\ref{polya}) are treated uniformly by the 
formula (\ref{eq_oneMinus}). For $\overline{v}\ne\overline{0}$, we 
could not do it and split (\ref{polya}) in formulas (\ref{eq_exactly1}) and 
(\ref{eq_odmoc}). A~selection of 
references discussing P\'olya's theorem is 
\cite{beck,bend_rich,bill,doyl_snel,fell,fost_good,
grim_wels,koch,lang,levi_pere,
nova,reny,reve,salv,wins,woes}. 

We obtain P\'olya's theorem (\ref{polya}) for 
$\overline{v}=\overline{0}$ via the symbolic method, that is, with the help of generating functions. The 
grid graph $G_d=\langle\Z^d,E_d\rangle$ is a~$2d$-regular transitive graph with the edges
$${\textstyle
\{\overline{a},\,\overline{b}\}\in E_d\iff
\sum_{i=1}^d|a_i-b_i|=1\,.
}
$$
We weight every edge $e\in E_d$ with $h(e)=\frac{1}{2d}$.
We extend the weight $h$ to the walks $w$ in $G_d$ and set $h(w)=(2d)^{-|w|}$ where $|w|$ ($\in\N_0=
\{0,1,\ds\}$) is the length of $w$, which is the number of 
edges in $w$. 
For $n\in\N_0$, we denote by $A_n$ the set of walks of length~$n$ in $G_d$ that start 
at the origin $\overline{0}$ and later 
revisit it. Then
$${\textstyle
P_d(\overline{0},\,n)=
\sum_{w\in A_n}h(w)=:h(A_n)
}
$$
and we introduce the main generating function
$${\textstyle
F_A(x)=\sum_{n\ge0}h(A_n)x^n\ \ (\in\C[[x]])\,.
}
$$
We also introduce three auxiliary generating functions 
$$
{\textstyle
F_B(x)=
\sum_{n\ge0}h(B_n)x^n,\ F_C(x)=\sum_{n\ge0}h(C_n)x^n,\ F_D(x)=\sum_{n\ge0}h(D_n)x^n\,.
}
$$
Here $B_n$ is the set of walks with length $n$ in $G_d$ that start and end 
at $\overline{0}$ (they are called closed walks), $C_n\sus B_n$ is the 
set of closed walks with inner vertices different from 
$\overline{0}$, and $D_n$ is the set of all walks of length $n$ in $G_d$ starting at~$\overline{0}$. We have $h(A_0)=h(C_0)=0$ and $h(B_0)=h(D_0)=1$. Clearly, $F_D(x)=
\frac{1}{1-x}$ because $h(D_n)=1$ for every $n\in\N_0$ (the weights 
of all length $n$ walks starting at $\overline{0}$ sum up to $1$).

The symbolic method yields formal identities 
$${\textstyle
F_A(x)=F_C(x)F_D(x)=F_C(x)\frac{1}{1-x}\,\text{ and }\,F_C(x)=1-\frac{1}{F_B(x)}\ \ (\text{in }\,\C[[x]])\,.
}
$$
We generalize them in Proposition~\ref{prop_12rel}. The first identity implies that  
$${\textstyle
h(A_n)=\sum_{j=0}^n h(C_j)\,. 
}
$$
Thus, the sequence $(h(A_n))_{n\ge0}$ is non-decreasing and therefore has 
a~limit; this is mentioned already in \cite{poly}.
We consider the sums
$${\textstyle
F_B(1):=\lim_{n\to\infty}\sum_{j=0}^n h(B_j)\,\text{ and }\,F_C(1):=\lim_{n\to\infty}\sum_{j=0}^n h(C_j)\ \ 
(\in\C)}\,.
$$
Since $0\le h(B_n),h(C_n)\le h(D_n)=1$ for every $n\in\N_0$, the formal power series 
$F_B(x)$ and $F_C(x)$ have coefficients in $[0,1]$. The second formal identity $F_C(x)=1-\frac{1}{F_B(x)}$ therefore implies the numeric identity
$${\textstyle
F_C(1)=1-\frac{1}{F_B(1)}\ \ (\in[0,\,1])\,,
}
$$
where $F_B(1)\in[1,\,+\infty)\cup\{+\infty\}$. We provide more details. Formally,  
$$
F_B(x)F_C(x)=F_B(x)-1\ \ (\text{in }\,\C[[x]])\,. 
$$
The power series $F_B(x)$ and 
$F_C(x)$ absolutely converge for every number 
$x$ in $[0,1)$. Using the operations of linear combination, grouping, and 
product of series, treated in respective Propositions~\ref{prop_LKseries}, 
\ref{prop_grouSer}, and \ref{prop_proSer}, we obtain for every number $x\in[0,1)$ the equality
$$
F_B(x)F_C(x)=F_B(x)-1\ \ (\text{in }\,\R_{\ge0})\,.
$$
Since $F_B(x)\ge1$ for every $x\in[0,1)$, we can divide by 
$F_B(x)$ and get the numeric 
identity 
$${\textstyle
F_C(x)=1-\frac{1}{F_B(x)}\,\text{ for every }\,x\in[0,\,1)\,.
}
$$
It follows that
$${\textstyle
F_C(1)=\lim_{x\nearrow1}F_C(x)=
\lim_{x\nearrow1}\big(1-\frac{1}{F_B(x)}\big)=
1-\frac{1}{\lim_{x\nearrow1}F_B(x)}=1-\frac{1}{F_B(1)}\,,
}
$$
where $\lim_{x\nearrow1}$ denotes that $x$ approaches $1$ via  
$x\in[0,1)$. The first equality follows from 
Theorem~\ref{thm_abel3} (Abel's theorem for nonnegative 
real coefficients). The second equality follows from the above 
numeric identity. The fourth and third equalities follow from 
Theorem~\ref{thm_abel3} and from the arithmetic of functional limits. Thus
$$
\lim_{n\to\infty}
P_d(\overline{0},\,n)=
\lim_{n\to\infty} h(A_n)=
\lim_{n\to\infty}{\textstyle\sum_{j=0}^n h(C_j)=F_C(1)=
1-\frac{1}{F_B(1)}\,.
}
$$
Since $F_B(1)=+\infty$ for $d\le2$ and $1<F_B(1)<+\infty$ for $d\ge3$,  
which we prove in Appendix~\ref{sec_sumsB}, P\'olya's theorem (\ref{polya}) 
for $\overline{v}=\overline{0}$ follows.\footnote{We 
include Appendix~\ref{sec_sumsB} to make our presentation of the generating functions proof of P\'olya's theorem for 
$\overline{v}=\overline{0}$ self-contained.}
\eproof

\noindent
See Theorem~\ref{thm_v1conNN} for a~generalization. The formula
\begin{equation}\label{eq_oneMinus}
{\textstyle
\lim_{n\to\infty}
P_d(\overline{0},\,n)=
\lim_{n\to\infty} h(A_n)=1-\frac{1}{F_B(1)}
}
\end{equation}
unifying the cases $d\le2$ and $d\ge3$ of the dichotomy (\ref{polya}) was mentioned in the 1990s by P.~Flajolet in an e-mail discussion group \cite{band_pers}. 

The use of generating functions in the proof of P\'olya's theorem is 
well known, but the theorem itself is usually viewed as a result on random 
walks in probability theory. See, for example, \cite{bill} for 
a~probabilistic proof. In our article, we demonstrate that 
P\'olya's theorem can also be viewed as a result in enumerative combinatorics. This 
interpretation does not require the technical tools of probability 
theory, such as Markov chains.

\subsection{P\'olya's theorem for 
$\overline{v}\ne\overline{0}$}\label{subsec_PTvSMvn0}

In this subsection, we proceed briefly. 
Let $n\in\N_0$ and  
$\overline{v}\in\Z^d\setminus\{\overline{0}\}$. We denote by $A_{\overline{v},n}$ the set of 
walks with length $n$ in $G_d$ starting at $\overline{0}$ and visiting 
$\overline{v}$. For $d\le2$ when $F_B(1)=+\infty$, the 
formula
\begin{equation}\label{eq_exactly1}
\lim_{n\to\infty}
P_d(\overline{v},\,n)=
\lim_{n\to\infty} h(A_{\overline{v},\,n})=1
\end{equation}
is a~particular case of part~1 of Theorem~\ref{thm_vne1conInf}.
For $d\ge3$ when $1<F_B(1)<+\infty$, a~more 
involved version of the argument in
the previous subsection produces the formula
\begin{equation}\label{eq_odmoc}
\lim_{n\to\infty}
P_d(\overline{v},\,n)=
\lim_{n\to\infty} h(A_{\overline{v},\,n})={\textstyle
\sqrt{1-\frac{F_{B,\,\overline{v}}(1)}{F_B(1)}}}\,.
\end{equation}
Here $F_{B,\overline{v}}(x)$ is the generating function of the walks in 
$G_d$ that begin and end at $\overline{0}$ and avoid $\overline{v}$. It is easy to see that
$$
1\le F_{B,\overline{v}}(1)<F_B(1)\,.
$$
This implies the second case of the dichotomy (\ref{polya}). We obtain a~more general result in 
Theorem~\ref{thm_vne1conNN}. We have not seen the formula (\ref{eq_odmoc}) in the literature.

\subsection{Our generalization of P\'olya's theorem}
\label{subsec_genPT}

In our opinion, P\'olya's theorem (\ref{polya}) divides into two parts. The first part, which we generalize 
in our article, consists of the numeric identities 
(\ref{eq_oneMinus})--(\ref{eq_odmoc}). The second part consists of 
the values $F_B(1)=+\infty$ for $d\le2$ and $F_B(1)\in(1,+\infty)$ 
for $d\ge3$. 
For the completeness of exposition, we present our proof of the second part in Appendix~\ref{sec_sumsB}. 

We generalize P\'olya's problem
as follows. We consider the countable complete directed graph 
$${\textstyle
\vec{K}_{\N}=\langle\N,\,\N_2\rangle:=\big\langle\N,\,
\N\times\N\setminus\{\langle n,\,n\rangle\cc\;n\in\N\}\big\rangle
}
$$
and endow its edges with weights $h\cc\N_2\to\C$ satisfying the convergence condition that for every $n\in\N_0$, the series 
$$
{\textstyle
\sum_{w\in D_n}h(w)
}
$$
absolutely converges. Here $D_n$ is the set of all walks of length $n$ in 
$\vec{K}_{\N}$ starting at the vertex $1$, and the weight $h$ is extended 
to walks in $\vec{K}_{\N}$ in the usual multiplicative way (stated 
formally 
before Definition~\ref{def_light}).
We say that a~weight $h\cc\N_2\to\C$ is {\em symmetric}, if $h(\langle 
u,v\rangle)=h(\langle v,u\rangle)$ for every $u,v\in\N$ with $u\ne v$.
For example, the constant weights in P\'olya's problem are symmetric. 

For $v\in\N$ and $n\in\N_0$, we denote by $A_{v,n}$ the set of walks of 
length $n$ in $\vec{K}_{\N}$ starting at~$1$ and later visiting the vertex 
$v$. It follows from the above convergence condition that the 
series $\sum_{w\in A_{v,n}}h(w)$ absolutely converges as well. We denote its 
sum by $h(A_{v,n})$, introduce the main generating function
$$
{\textstyle
F_{A,\,v,\,h}(x)=\sum_{n\ge0}h(A_{v,n})x^n\,,
}
$$
and investigate the behavior of $F_{A,v,h}(x)$ and its coefficients.
This is a~far-reaching generalization 
not only of P\'olya's problem, but of every Markov chain on $\N$ 
starting at~$1$.

The novel features are: (i) the series like $\sum_{w\in 
D_n}h(w)$ and $\sum_{w\in A_{v,n}}h(w)$ are for $n>1$ infinite 
and (ii) the weights $h$ have values generally outside of $[0,1]$, which 
means that we transcend probabilistic weights. We generalize 
the numeric identities (\ref{eq_oneMinus})--(\ref{eq_odmoc}) 
to the sums $h(A_{v,n})$ and state these generalizations in terms of 
four quantities listed below in 
the order of decreasing precision. 
\begin{enumerate}
\item The limits $\lim_{n\to\infty}h(A_{v,n})$.  
\item The sums $F_{A,v,h}(1)_a$ of the absolutely convergent series 
$${\textstyle
\sum_{n\ge0}h(A_{v,\,n})=
h(A_{v,0})+h(A_{v,1})+\ds\,.
}
$$
\item The sums $F_{A,v,h}(1)_c$ of the same series if it converges only 
conditionally. 
\item The generalized sums 
$$
\lim_{x\nearrow1}F_{A,\,v,\,h}(x)=
\lim_{x\nearrow1}{\textstyle
\sum_{n\ge0}h(A_{v,\,n})x^n
}\,.
$$
\end{enumerate} 
We twist these quantities by powers of a~number $\xi\in\C$: if 
$${\textstyle
U(x)=\sum_{n\ge0}u_nx^n
}
$$ 
is a~complex formal power series, we 
define its $\xi$-twist by
$${\textstyle
(U)_{\xi}(x):=U(\xi x)=
\sum_{n\ge0}\xi^nu_nx^n\,. 
}
$$
There are 
two reasons for working with $\xi$-twists. (i) In the probabilistic 
approach, the weights of all edges incident to a~fixed vertex sum up 
to $1$. We generalize this condition either by dropping it or by 
requiring that the sums are equal to $\xi$. (ii)
We generalize the weight ranges 
$\R_{\ge0}$ and $\R$ to the complex half-lines $\xi\R_{\ge0}$ 
and lines $\xi\R$. We call such weights {\em 
complex nonnegative} and {\em complex real}, respectively. Both (i) and (ii) lead to $\xi$-twists. In   
Sections~\ref{sec_extTwoPol} and~\ref{sec_extTwoPol2} we obtain 
along these lines a total of sixteen 
general P\'olya's theorems. In the probabilistic approach, when weights 
are real non-negative, the absolute and conditional sums of the series 
mentioned in items~2 and~3 coincide. 

\subsection{Plan of the article and perspectives}\label{subsec_plan}

Let $h\cc\N_2\to\C$. 
We extend only a fragment of the symbolic method in enumerative 
combinatorics. In Proposition~\ref{prop_12rel}, we give 
detailed proofs for only two formal 
identities
$${\textstyle
F_{A,\,h}(x)=F_{C,\,h}(x)F_{D,\,h}(x)\,\text{ and }\, 
F_{C,\,h}(x)=1-\frac{1}{F_{B,\,h}(x)}}\ \ (\text{in }\C[[x]])\,.
$$
Here $F_{A,\,h}(x):=F_{A,\,1,\,h}(x)$
and so on. They look similar to the two formal identities used in Subsection~\ref{subsec_PTvSMv0}, 
but they are, in fact, new and {\em need to be proven anew}. Coefficients of the generating 
functions involved in them arise as sums of infinite 
series. The coefficients in Subsection~\ref{subsec_PTvSMv0} are finite sums or 
finite linear combinations.

So we need a~theory of series defined on arbitrary countable sets, 
like the series $\sum_{w\in A_{v,n}}h(w)$. We could not find such a~theory in the literature in 
a~convenient form, 
and we develop it in Section~\ref{sec_serGF}. For example, 
we lift the ``obvious'' product identity $|A\times B|=|A|\cdot|B|$ for finite 
sets $A$ and $B$, which is frequently used in the symbolic method, to
Proposition~\ref{prop_proSer} on sums of products of absolutely convergent 
series. Proposition~\ref{prop_grouSer} 
concerns another important operation on series, the grouping. Other 
results of interest in Section~\ref{sec_serGF} are 
three Abel's Theorems~\ref{thm_abel}--\ref{thm_abel3} (stated without 
proofs) and several conversions of formal identities to numeric 
ones, like Proposition~\ref{prop_bezAbel3}.

In Section~\ref{sec_extTwoPol} we obtain general 
P\'olya's theorems for $v=1$, that is, when the vertex to be (re)visited 
by walks is the starting one. We obtain seven general P\'olya's 
theorems in this case: \underline{Theorems~\ref{thm_v1gen1}--\ref{thm_v1conInf}}. Part~$4$ of 
Theorem~\ref{thm_v1conNN} is the classical criterion of transience of 
walks with probabilistic ($1$-convex) weights $h\cc\N_2\to\R_{\ge0}$. Theorems~\ref{thm_v1con1}--\ref{thm_v1conInf} generalize this criterion to complex 
weights.
In Section~\ref{sec_extTwoPol2},  
we obtain analogous 
\underline{Theorems~\ref{thm_vne1gen1}--\ref{thm_vne1conInf}}. Now $v\ne1$, 
and the theorems are more complicated. In 
Theorems~\ref{thm_v1genInf}, \ref{thm_v1conNN}, 
\ref{thm_v1conInf}, \ref{thm_vne1genInf}, 
\ref{thm_vne1conNN}, and  \ref{thm_vne1conInf} where weights are nonnegative 
or complex nonnegative and the subseries bound is therefore available (see 
Subsection~\ref{subsec_series}), we allow any vertex $u\in\N$ to be the 
starting one for walks. 
The proofs of the three formal identities 
\begin{enumerate}
\item $F_{A,v,h}(x)
=F_{C,v,h}(x)F_{D,h}(x)$,
\item $F_{B,\,h}(x)=F_{B,v,h}+F_{C,v,h}(x)^2F_{B,h}(x)$, and
\item $F_{C,v,h}(x)=F_{B,v,h}(x)
F_{E,v,h}(x)$
\end{enumerate}
in Proposition~\ref{prop_12relvne1} 
are only sketched. If $v\ne1$, one can introduce more auxiliary 
generating functions than
in the case $v=1$. In Proposition~\ref{prop_onFv}, we 
present two additional formal identities for the generating function 
$F_{F,v,h}(x)$ of walks going from $1$ to $v$. In 
Corollary~\ref{cor_monoton} we utilize the square in identity~2 and show that for every weight $h\cc\N_2\to\R$, every 
vertex $v\ne1$, every number  $x\in(-1,1)$, and with appropriate assumptions, the following holds. If $F_{B,h}(x)=0$, then $F_{B,v,h}(x)=0$. If $F_{B,h}(x)\ne0$, then
$$
\frac{F_{B,\,v,\,h}(x)}{F_{B,\,h}(x)}\le1\,.
$$
This generalizes the obvious bound 
$$
0\le F_{B,\,v,\,h}(x)\le F_{B,\,h}(x)\ \  (x\in[0,1))
$$ 
for nonnegative weights $h\cc\N_2\to\R_{\ge0}$. Recall that $F_{B,h}(x)$ is the generating function of weighted closed walks and that 
$F_{B,v,h}(x)$ counts weighted closed walks avoiding~$v$.  
Section~\ref{sec_concl} contains concluding remarks of various kinds.

We hope to continue our work in four directions. First, it would be nice 
to really have an extension of the symbolic method in enumerative 
combinatorics to infinite coefficient series \cite{klaz0}. However, to construct it in a general form 
one needs to gain experience from 
concrete applications, and the present article is only a beginning. Second, 
in \cite{klaz_hors2} we replace complex weights with weights taken from a field endowed with 
a complete non-Archimedean norm. Third, one can view the formal identities in the present article, 
such as the above $F_{A,h}(x)=
F_{C,h}(x)F_{D,h}(x)$, as living in a~semiformal algebra of formal power 
series, an algebra that allows coefficients obtained via limit 
transitions. We hope to develop this approach to 
formal power series in \cite{klaz1}. Fourth, as we already 
mentioned, our approach generalizes not only the rather particular P\'olya's 
problem but even every Markov chain on $\N$ starting at $1$. In our 
generalization to edge-weighted $\vec{K}_{\N}$, we avoid the elephant 
in the room of Markov chains, which is the underlying uncountable (!) 
probability space. We look into this matter in 
\cite{klaz2}.   

\section{Complex series and generating functions}\label{sec_serGF}

In this section, we develop a theory of complex series defined 
on at most countable sets and obtain conversions of formal identities 
satisfied by generating functions to numeric identities at $x=1$. These 
results are needed in Sections~\ref{sec_extTwoPol} and~\ref{sec_extTwoPol2}. We believe that some of our results on both topics 
are original contributions. For example, we are not aware of any prior satisfactory treatment of the grouping operation for series in the literature; see 
Proposition~\ref{prop_grouSer} and the preceding definitions here.

\subsection{Series and light weights}\label{subsec_series}

A~{\em series}  is a~map 
$h\cc X\to\C$ defined on a~finite or countable set $X$. 
We also write it as 
$\sum_{x\in X}h(x)$. We say that
{\em $\sum_{x\in X}h(x)$ absolutely converges} if two equivalent conditions hold.

\begin{enumerate}
\item For some constant $c>0$ and every finite set $Y\sus X$, the bound 
$${\textstyle
\sum_{x\in Y}|h(x)|\le c
}
$$
holds.
\item For every bijection $f\cc\N\to 
X$, the limit 
$$
s=\lim_{n\to\infty}{\textstyle\sum_{i=1}^n h(f(i)) \ \ (\in\C)
}
$$
exists and does not depend on $f$.
\end{enumerate}
For finite sets $X$, condition~2 is vacuously true and finite series absolutely converge.  
The independence of the limit $s$ on the bijection $f$ may be omitted from condition~2 because if 
$f_1,f_2\cc\N\to X$ are bijections such that 
$${\textstyle
\lim_{n\to\infty}\sum_{j=1}^n h(f_i(j))=s_i\,\text{ and }\,s_1\ne s_2\,,
}
$$
we can construct a~bijection $f_3\cc\N\to 
X$ such that $\lim_{n\to\infty}\sum_{j=1}^n h(f_3(j))$ does not exist. 
We leave the proof of the equivalence of conditions~1 and~2 as an exercise for the interested reader.
We call the number $s$ the {\em sum} of the series $\sum_{x\in X}h(x)$ and denote it again by 
$\sum_{x\in X}h(x)$ or by $h(X)$ (we already used this notation). If the 
set $X=\{x_1,x_2,\ds,x_n\}$
is finite and nonempty, then every series  
$\sum_{x\in X}h(x)$ absolutely converges and has 
the sum 
$$
h(x_1)+h(x_2)+\ds+h(x_n)\,. 
$$
For $n=0$ we set $h(\emptyset)=0$. If 
$${\textstyle
U=U(x)=\sum_{n\ge0}u_nx^n\ \ (\in\C[[x]])
}
$$ 
is a~formal power series, so that $x$ is 
a~formal variable, we denote for any $c\in\C$ by $U(c)$ the 
sum of the power series
$\sum_{n\ge0}u_nc^n$,
if it absolutely converges.
Let
$$
H_1:={\textstyle
\{U(x)\in\C[[x]]\cc\;\text{$U(c)$ absolutely converges for every $c\in\C_{<1}$}\}\,,
}
$$
where $\C_{<1}:=\{c\in\C\cc\;|c|<1\}$. So for any formal power series 
$U(x)\in H_1$, the sums $U(c)$ define a~function $U\cc\C_{<1}\to\C$ by 
$c\mapsto U(c)$. 

Let $S=\sum_{x\in X}h(x)$ be a~series. If $Y\sus X$, we say that $\sum_{x\in 
Y}h(x)$ is a~{\em subseries} of $S$. Any subseries of an 
absolutely convergent series absolutely converges. In general, we cannot compare the sums of a~series and its 
subseries: it may very well happen that $h(X)=0$, but $h(Y)\ne0$. However, if $h\cc X\to\R_{\ge0}$, 
$\sum_{x\in X}h(x)$ 
absolutely converges, and if $Y\sus X$, then the {\em subseries bound}
$$
0\le h(Y)\le h(X)
$$
holds.

We formally introduce walks in $\vec{K}_{\N}$ and extensions of weights on $\N_2$ to walks. A~{\em walk $w$} in $\vec{K}_{\N}$ of {\em 
length $|w|=n$},
with $n\in\N_0$, is any $(n+1)$-tuple
$$
w=\langle u_0,\,u_1,\,\ds,\,u_n\rangle\in\N^n
$$
of vertices in $\vec{K}_{\N}$
satisfying $u_{i-1}\ne u_i$, so that 
$\langle u_{i-1},u_i\rangle\in\N_2$ for every $i=1,\ds,n$. Recall that
$\N_2=\{\langle m,n\rangle\cc\;m,n\in\N, m\ne n\}$.
If $h\cc\N_2\to\C$ is a~{\em weight}, we define $h(w):=1$ if $|w|=0$, and 
$$
{\textstyle
h(w):=\prod_{i=1}^n h(\langle u_{i-1},u_i\rangle)\,\text{ if }\,|w|=n>0\,.
}
$$
Recall that $D_n$ is the 
set of walks $w$ of length~$n$ in $\vec{K}_{\N}$ starting at $u_0=1$.

\begin{defi}[light weights]\label{def_light}
A~weight $h\cc\N_2\to\C$ is light if, for every $n$ in $\N_0$, the series
$$
{\textstyle
\sum_{w\in D_n}h(w)
}
$$
absolutely converges.    
\end{defi}

We show that the starting vertex $1$ may be replaced with any other 
starting vertex $u$ that is reachable from $1$ by a~walk with nonzero 
weight.

\begin{prop}\label{prop_moveStVer}
Let $h\cc\N_2\to\C$ be a~light weight, let $u\in\N$ be such that 
there exists a~walk $w$ in $\vec{K}_{\N}$ from $1$ to $u$ with 
$h(w)\ne0$, and for $n\in\N_0$ let $D_{u,n}$ be the set of walks of 
length $n$ in 
$\vec{K}_{\N}$ starting 
at $u$. Then for every $n\in\N_0$, the series
$$
{\textstyle
\sum_{w\in D_{u,\,n}}h(w)
}
$$
absolutely converges.
\end{prop}
\proof
Let $w_0$ be a~walk from $1$ to $u$ with nonzero weight, and let $m=|w_0|$ ($\in\N_0$) and 
$c=h(w_0)\ne0$. Suppose that $h$ is light and that $n\in\N_0$. For every 
$w\in D_{n,u}$ we define a~concatenated walk $w'=w_0w\in D_{m+n}$. Then 
the equality of series 
$$
{\textstyle
\sum_{w\in D_{n,\,u}}h(w)=
c^{-1}\sum_{w\in D_{n,\,u}}h(w')
}
$$
holds. The latter series is a~subseries of the absolutely 
convergent series 
$${\textstyle
\sum_{w\in D_{m+n}}h(w)\,. 
}
$$
Therefore, both the 
latter and the former series absolutely converge. We use the linear combination of series treated next in Proposition~\ref{prop_LKseries}. 
\eproof

\subsection{Operations with series}

Let $R=\sum_{x\in X}g(x)$ and $S=\sum_{x\in X}h(x)$ be 
series and let $\al,\be\in\C$. The series
$$
{\textstyle
\al R+\be S=\sum_{x\in X}(\al g(x)+\be h(x))
}
$$
is the {\em linear combination} of $R$ and $S$. We omit the 
easy proof of the next result.

\begin{prop}\label{prop_LKseries}
If $R$ and $S$ absolutely converge and have sums $r$ and $s$, 
respectively, then $\al R+\be S$ absolutely
converges and has sum $\al r+\be s$.
\end{prop}
This proposition easily extends to linear combinations with more than 
two terms. The next approximation lemma is useful. 

\begin{lem}\label{lem_apprLem}
If $S=\sum_{x\in X}h(x)$ absolutely converges and has sum $s$, then for every $\ep>0$ there is a~finite set 
$Y\sus X$, denoted by 
$Y(S,\,\ep)$, such that for every finite set $Z$ with $Y\sus Z\sus X$,  the bound 
$$
|s-h(Z)|\le\ep
$$
holds.
\end{lem}
\proof
Let $\ep>0$ be given.
For finite $X$ we set $Y=X$. For countable $X$ we take any bijection
$f\cc\N\to X$, take $N\in\N$ such that 
$${\textstyle
\big|\sum_{n=1}^N h(f(n))-
s\big|\le\frac{\ep}{2}\,\text{ and }\,\sum_{n>N}|h(f(n))|\le
\frac{\ep}{2} \,,
}
$$
and set $Y=\{f(1),f(2),\ds,f(N)\}$. Then for every finite set $Z$ with $Y\sus Z\sus X$, indeed 
$$
{\textstyle
|s-h(Z)|\le
|s-h(Y)|+
\sum_{x\in Z\setminus Y}|h(x)|\le
\frac{\ep}{2}+\frac{\ep}{2}=\ep
\,.
}
$$
\eproof

Two important operations with series are grouping and product. A~{\em 
partition} of a~set $X$ is a~set $P$ of nonempty and
mutually disjoint sets such that $\bigcup P=X$. If $X$ is at most countable, 
then so is $P$ and every set $Z\in P$. If $S=\sum_{x\in X}h(x)$ is 
a~series and $P$ is a~partition of $X$ such that for 
every set $Z\in P$ the subseries $S_Z=\sum_{x\in Z}h(x)$ absolutely 
converges and has 
sum $s_Z$ ($\in\C$), then the series
$$
%{\textstyle
S_P=\sum_{Z\in P}s_Z
%}
$$
is called the $P$-{\em grouping}, or just {\em grouping}, of $S$. 

\begin{prop}\label{prop_grouSer}
Suppose that $S=\sum_{x\in X}h(x)$ absolutely
converges and has sum $s$, and that $P$ is a~partition of $X$. Then the 
$P$-grouping of $S$ is correctly defined, absolutely converges, and 
has the same sum $s$. 
\end{prop}
\proof
For every $Z\in P$, the series
$S_Z=\sum_{x\in Z}h(x)$ is a~subseries of $S$ and absolutely 
converges. So the sum $s_Z$ exists and the series $S_P$ is correctly defined. We show that it 
absolutely converges. Let $c>0$ be a~constant witnessing the absolute 
convergence of $S$ and let 
$\{Z_1,Z_2,\ds,Z_n\}\sus P$ be a~finite set.
Using Lemma~\ref{lem_apprLem}, we take for $i=1,2,\ds,n$ finite sets 
$$
Z_i'=Y(S_{Z_i},\,2^{-i})\ \  (\sus Z_i) 
\,\text{ and }\,Z_0=Z_1'\cup Z_2'\cup\ds\cup Z_n'\ \ (\sus X)\,.
$$
Then
\begin{eqnarray*}
{\textstyle
\sum_{i=1}^n|h(Z_i)|
}&\le&   
{\textstyle
\sum_{i=1}^n|h(Z_i)-h(Z_i')|+
\sum_{x\in Z_0}|h(x)|}\\
&\le&{\textstyle\sum_{i=1}^n
2^{-i}+c\le 1+c}
\end{eqnarray*}
and we see that the series $S_P$ absolutely converges. 

Let $t$ be the sum of $S_P$.
It suffices to show that $|t-s|\le\ep$ for every $\ep>0$. 
Let $\ep>0$ be given. We use 
Lemma~\ref{lem_apprLem} and take finite sets 
$${\textstyle
X'=Y(S,\,\frac{\ep}{3})\ \  (\sus X)
\,\text{ and }\,P'=Y(S_P,\,\frac{\ep}{3})\ \  (\sus P)\,. 
}
$$
We take a~finite set $\{Z_1,Z_2,\ds,Z_n\}\sus P$ such that 
$${\textstyle
P'\sus\{Z_1,\,Z_2,\,\ds,\,Z_n\}\,\text{ and }\,X'\sus\bigcup_{i=1}^n Z_i\,.
}
$$
We again use Lemma~\ref{lem_apprLem} and take for $i=1,2,\ds,n$ finite sets 
$${\textstyle
Z_i'=Y(S_{Z_i},\,2^{-i}\cdot\frac{\ep}{3})\ \  (\sus Z_i)\,. 
}
$$
For $i=1,2,\ds,n$ we define 
$$
Z_i''= Z_i'\cup(X'\cap Z_i)\ \ (\sus Z_i)
$$
and form the disjoint union 
$${\textstyle
X_0=\bigcup_{i=1}^n Z_i''\ \  
(\sus X)\,. 
}
$$
Then $X_0$ is a~finite set and $X'\sus X_0$. Also, $Z_i'\sus Z_i''$ for every $i=1,2,\ds,n$. It follows that
\begin{eqnarray*}
|s-t|&\le&{\textstyle|s-h(X_0)|+
\sum_{i=1}^n|h(Z_i'')-h(Z_i)|+|\sum_{i=1}^n h(Z_i)-t|}\\
&\le&{\textstyle\frac{\ep}{3}+\sum_{i=1}^n 2^{-i}\cdot\frac{\ep}{3}+\frac{\ep}{3}\le\ep\,.
}
\end{eqnarray*}
\eproof

If $S_i=\sum_{x\in X_i}h_i(x)$ for $i=1,2$ are two series, then the series 
$${\textstyle
S_1\otimes S_2:=\sum_{\langle x,\,y\rangle\in 
X_1\times X_2}h_1(x)\cdot h_2(y)
}
$$ 
is the {\em product} of $S_1$ and $S_2$. 

\begin{prop}\label{prop_proSer}
If $S_1$ and $S_2$ absolutely converge and have sums $s_1$ and 
$s_2$, respectively, then the product series $S_1\otimes S_2$ 
absolutely converges and has the sum $s_1s_2$.
\end{prop}
\proof
We show that $S_1\otimes S_2$ absolutely converges. Let 
$c>0$ be such that for any finite sets $Y_i\sus 
X_i$ and $i=1,2$ we have
$|h_i|(Y_i)\le c$. 
Let $A\sus X_1\times X_2$ be a~finite set. We 
take finite sets $Y_i\sus X_i$ such that $A\sus Y_1\times Y_2$ and get the desired bound:
$${\textstyle
\sum_{\langle x,\,y\rangle\in A}|h_1(x)h_2(y)|\le
\sum_{x\in Y_1}|h_1(x)|\cdot
\sum_{y\in Y_2}|h_2(y)|\le c\cdot c=c^2\,.
}
$$

We show that the sum $s$ of  $S_1\otimes S_2$ equals $s_1s_2$. 
It suffices to prove that for every $\ep\in(0,1)$ we have $|s-
s_1s_2|\le\ep$. So let $\ep\in(0,1)$. 
Using Lemma~\ref{lem_apprLem} we take finite sets 
$${\textstyle
X'_1=Y\big(S_1,\,\frac{\ep}{4(1+|s_2|)}\big)\ \  (\sus X_1),\ \  
X'_2=Y\big(S_2,\,\frac{\ep}{4(1+|s_1|)}\big)\ \  (\sus X_2)\ 
}
$$
and 
$${\textstyle
Z=Y(S_1\otimes 
S_2,\frac{\ep}{4})\ \ (\sus X_1\times X_2)\,. 
}
$$
For $i=1,2$ we take finite sets
$X_i''$  such that $X_i'\sus X_i''\sus X_i$ and $Z\sus X_1''\times X_2''$. Then 
\begin{eqnarray*}
|s-s_1s_2|&\le&
{\textstyle
\big|s-\sum_{\langle x,\,y\rangle\in X_1''\times X_2''}h_1(x)h_2(y)\big|+
|h_1(X_1'')h_2(X_2'')-s_1s_2|}\\
&\le&
{\textstyle\frac{\ep}{4}+|(s_1+t_1)(s_2+t_2)-s_1s_2|\,,
}
\end{eqnarray*}
where 
$|t_1|\le\frac{\ep}{4(1+|s_2|)}$ and $|t_2|\le\frac{\ep}{4(1+|s_1|)}$. Hence $|s-s_1s_2|\le
\frac{\ep}{4}+|s_1t_2|+|s_2t_1|+|t_1t_2|
\le\frac{\ep}{4}+\frac{\ep}{4}+
\frac{\ep}{4}+\frac{\ep}{4}=\ep$.
\eproof

\noindent
Let $S_i=\sum_{x\in X_i}h_i(x)$ for  
$i=1,2,\ds,k$ be $k$ series. The {\em product} of $S_1$, $S_2$, $\ds$, $S_k$ is the series
$$
{\textstyle
S_1\otimes S_2\otimes\ds\otimes S_k=\sum_{\overline{x}\in\overline{X}}h_1(x_1)\cdot h_2(x_2)\cdot\ldots\cdot h_k(x_k)\,,
}
$$ 
where $\overline{x}=\langle x_1,x_2,\ds,x_k\rangle$ and 
$\overline{X}=X_1\times X_2\times\ds\times X_k$. By iterating 
Proposition~\ref{prop_proSer} we obtain the following result.

\begin{prop}\label{prop_prodKser}
Let $S_1$, $S_2$, $\ds$, $S_k$ be absolutely convergent series with respective sums $s_1$, $s_2$, $\ds$, $s_k$. Then the series 
$S_1\otimes S_2\otimes\ds\otimes S_k$ absolutely converges and has sum $s_1s_2\ds s_k$.
\end{prop}
Interestingly, a~partial reversal of Proposition~\ref{prop_proSer} holds.

\begin{cor}\label{cor_InvproSer}
Let $S_1$ and $S_2$ be nonempty series whose summands are not identically zero. If $S_1\otimes S_2$ absolutely 
converges and has sum $s$, then $S_1$ and $S_2$ 
absolutely converge, and their respective sums $s_1$ and $s_2$ satisfy $s_1s_2=s$.
\end{cor}
\proof
This follows from Proposition~\ref{prop_proSer} if we 
show that $S_1$ and $S_2$ absolutely converge. Let $S_i=\sum_{x\in 
X_i}h_i(x)$ for $i=1,2$ and let $c>0$ be such that for every finite set $X\sus X_1\times X_2$ we have 
$$
{\textstyle\sum_{\langle x,\,y\rangle\in X}
|h_1(x)h_2(y)|\le c\,.
}
$$
We show that $S_1$  absolutely 
converges. We take any $y_0\in X_2$ with $h_2(y_0)\ne0$. Then for 
every finite set $X\sus X_1$, the bound 
$$
{\textstyle
|h_1|(X)=\frac{1}{|h_2(y_0)|}\sum_{\langle x,\,y\rangle\in X\times\{y_0\}}|h_1(x)h_2(y)|
\le\frac{c}{|h_2(y_0)|}
}
$$
holds.
Thus $S_1$ absolutely converges. Similarly for $S_2$.
\eproof

\noindent
It is clear that the assumptions on the series $S_1$ and $S_2$ cannot be 
omitted. 

\subsection{Slim weights}

Let $h\cc\N_2\to\C$ be a~weight.
We obtain a sufficient condition for the lightness of $h$ in terms of the 
values $h(e)$. 

\begin{defi}[sets $V(h,n)$ and $V(h)$]
Let $h\cc\N_2\to\C$ be a~weight. For $n$ in $\N_0$, we denote by $V(h,n)$ the 
set of vertices $u\in\N$ such that there exists a~walk $w$ in $\vec{K}_{\N}$ 
from $1$ to $u$ such that $|w|\le n$ and $h(w)\ne0$. We set 
$${\textstyle
V(h)=\bigcup_{n\ge0}V(h,n)\,.
}
$$
\end{defi}

\begin{defi}[slim weights]\label{def_slimWei}
A~weight $h\cc\N_2\to\C$ is slim if for every $n$ 
in $\N_0$ there is 
a~constant $c(n)>0$ such that for every vertex 
$u\in V(h,n)$ and every finite set 
$X\sus\N\setminus\{u\}$, the bound
$${\textstyle
\sum_{v\in X}|h(\langle u,\,v\rangle)|\le c(n)
}
$$
holds.
\end{defi}

\begin{prop}\label{prop_slimWei}
Every slim weight  $h\cc\N_2\to\C$ is light.
\end{prop}
\proof
Let $h\cc\N_2\to\C$ be slim. We show by induction on $n\in\N_0$ that the series 
$\sum_{w\in D_n}h(w)$ 
absolutely converges. For $n=0$ it is trivial because $D_0=\{\langle1\rangle\}$. 
Let $n>0$ and let
$c>0$ be a~constant such that 
$${\textstyle
\sum_{w\in X}|h(w)|\le c\,\text{ and }\,
\sum_{v\in X'}|h(\langle u,\,v\rangle)|\le c
}
$$ 
for every finite set $X\sus D_{n-1}$, every vertex $u\in V(h,n-1)$, and every finite set $X'\sus\N\setminus\{u\}$. 

Let $X\sus D_n$  
be a~finite set. We may assume that 
all walks in $X$ have nonzero weights. We decompose 
every walk $w\in X$ as $w=w'\,\ell(w)$, where $w'\in D_{n-
1}$, and $\ell(w)$ is the last vertex of $w$. Let $Y$ be the set of the walks $w'$. For $w'\in Y$, let 
$$
X(w')=\{v\in\N\setminus
\{\ell(w')\}\cc\;w'v\in X\}\,.
$$
We get the desired bound
\begin{eqnarray*}
{\textstyle
\sum_{w\in X}|h(w)|}&=&{\textstyle
\sum_{w'\in Y}|h(w')|\cdot\sum_{v\in X(w')}|h(\langle\ell(w'),\,v\rangle)|
}\\
&\le&{\textstyle
\sum_{w'\in Y}|h(w')|\cdot c\le c\cdot c=c^2\,.
}
\end{eqnarray*}
\eproof

\subsection{Convex weights}

To define $\xi$-convex weights, we need 
a~lemma.

\begin{lem}\label{lem_convFromLigh}
Let $h\cc\N_2\to\C$ be a~light weight and let $u\in V(h)$. Then the series 
$${\textstyle
\sum_{v\in\N\setminus\{u\}}
h(\langle u,v\rangle)
}
$$ 
absolutely converges.   
\end{lem}
\proof
There is a~nonzero walk $w_0=\langle 
u_0,u_1,\ds,u_n\rangle$ such that $u_0=1$ and $u_n=u$. Let $c=h(w_0)$ ($\ne0$) and let $W$ be the set of walks 
$\langle u_0,\ds,u_n,u_{n+1}\rangle$  in $\vec{K}_{\N}$.
So the vertices $u_0$, $u_1$, $\ds$, $u_n$ are fixed, and $u_{n+1}$ ranges in 
$\N\setminus\{u_n\}=\N\setminus\{u\}$. The series $\sum_{w\in 
W}h(w)$ is a~subseries of $\sum_{w\in D_{n+1}}h(w)$ 
and therefore absolutely 
converges. By 
Proposition~\ref{prop_LKseries}, the series
$$
{\textstyle
\sum_{v\in\N\setminus\{u\}}
h(\langle u,v\rangle)=\frac{1}{c}\sum_{w\in W}h(w)
}
$$
absolutely converges. 
\eproof

\begin{defi}[$\xi$-convexity]\label{def_conv}
Let $\xi\in\C$. We call a~light weight $h\cc\N_2\to\C$  
$\xi$-convex if 
for every vertex $u\in V(h)$, the sum
$${\textstyle
\sum_{v\in\N\setminus\{u\}}
h(\langle u,\,v\rangle)\,,
}
$$
which exists by Lemma~\ref{lem_convFromLigh}, equals $\xi$.
\end{defi}
We allow $\xi=0$ and consider $0$-convex light weights; they are 
related to the networks satisfying Kirchhoff's second law. Since $V(h)$ 
is the set of vertices that can be reached from $1$ by a~nonzero walk, the $\xi$-convexity of $h$ is not generally preserved if we replace $1$
with another starting vertex $u\in\N\setminus\{1\}$. However, it is clearly preserved if $u\in V(h)\setminus\{1\}$.

The second half of Section~\ref{sec_serGF}
deals with generating functions. 

\subsection{Two kinds of sums of generating functions at~$x=1$}\label{subsec_u11u12}

We use  univariate formal 
power series 
$${\textstyle
U=U(x)=\sum_{n\ge0}u_nx^n 
}
$$ 
with complex coefficients $u_n$. They form an integral domain $\C[[x]]$. 
For $U=\sum_{n\ge0}u_nx^n$ and $V=\sum_{n\ge0}v_nx^n$ in $\C[[x]]$, 
$$
{\textstyle
U+V=\sum_{n\ge0}(u_n+v_n)x^n\,\text{ and }\, U\cdot V=UV=\sum_{n\ge0}\big(\sum_{j=0}^n u_jv_{n-j}\big)x^n\,.
}
$$
Neutral elements are 
$0=0x^0+0x^1+\ds$ and $1=1x^0+0x^1+0x^2+\cdots$. 
As for additive inverses,
$${\textstyle
-U=\sum_{n\ge0}(-u_n)x^n\,.
}
$$
Units in 
$\C[[x]]$ are exactly the formal power series 
$\sum_{n\ge0}u_nx^n$ with $u_0\ne0$. 

Let $\C_1:=\{x\in\C\cc\;|x|=1\}$
and $\xi\in\C_1$. Recall that
$$
\xi\R_{\ge0}=\{\xi x\cc\;
x\in\R_{\ge0}\}\ \ (\sus\C)
$$
is a~complex ray going from the origin in 
the direction $\xi$. This set is closed under addition. If 
$\al_j\in\xi_j\R_{\ge0}$ for $j=1,2$, then $\al_1\al_2\in
(\xi_1\xi_2)\R_{\ge0}$. If
$$
h\cc\N_2\to\xi\R_{\ge0}
$$
is a~light weight, $n\in\N_0$, and $X\sus D_n$,  then $h(X)\in 
\xi^n\R_{\ge0}$. The complex lines 
$\xi\R$ go through the origin in the direction $\xi$ and have 
the same properties. Recall that for $\xi\in\C$ and 
$U=U(x)=\sum_{n\ge0}u_nx^n\in\C[[x]]$, the {\em $\xi$-twist} of $U(x)$ is 
$${\textstyle
(U)_{\xi}=(U)_{\xi}(x)=U(\xi x)=
\sum_{n\ge0}u_n\xi^n x^n\ \ (\in\C[[x]])\,.
}
$$
Note that $(U+V)_{\xi}=(U)_{\xi}+(V)_{\xi}$ and 
$(U\cdot V)_{\xi}=(U)_{\xi}\cdot (V)_{\xi}$. 

\begin{defi}[$U(1)_c$ and $U(1)_a$]\label{def_U12}
Let $U(x)=\sum_{n\ge0}u_nx^n$ 
be in $\C[[x]]$. We introduce the following terminology.
\begin{enumerate}
\item If the limit
$$
U(1)_c:=\lim_{n\to\infty}{\textstyle
\sum_{j=0}^nu_j\in\C
}
$$
exists, we say that the conditional sum $U(1)_c$ exists.
\item If the series 
$${\textstyle
U(1)=\sum_{n\ge0}u_n1^n=\sum_{n\ge0}u_n=u_0+u_1+\ds
}
$$ 
absolutely converges, we denote its sum by $U(1)_a$ and say that the 
absolute sum $U(1)_a$ ($\in\C$) exists. 
\end{enumerate}
\end{defi}
If $U(1)_a$ exists, then $U(1)_c$ exists and $U(1)_c=U(1)_a$. The 
opposite is, in general, not true. A~well known example is
$$
{\textstyle
U(x)=\sum_{n\ge1}(-1)^{n-1}\frac{1}{n}x^n}\,,
$$
when $U(1)_c=
\log 2$, but $U(1)_a$ does not exist. But if
$$
{\textstyle
U(x)=\sum_{n\ge0}u_nx^n\in\R_{\ge0}[[x]]
}
$$
and $U(1)_c$ exists, 
then $U(1)_a$ exists and $U(1)_a=U(1)_c$. The next proposition 
follows from Hadamard's formula for the radius of convergence of power 
series. Recall that $H_1$ is the set of $U(x)\in\C[[x]]$ such that the 
power series $U(c)$ absolutely converges for every number $c\in\C$ with $|c|<1$.

\begin{prop}\label{prop_ConpowSer}
Let $U(x)=\sum_{n\ge0}u_nx^n
\in\C[[x]]$. Then $U(x)\in H_1$ if and only if 
$$
\limsup_{n\to\infty}
|u_n|^{1/n}\le1\,. 
$$
In particular, existence of the conditional sum $U(1)_c$ is 
a~sufficient (but not 
necessary) condition for $U(x)\in H_1$.
\end{prop}

\noindent
Note that if $U\in\C[[x]]$ then
$(U)_{\xi}\in H_1$ either for all or no $\xi\in\C_1$. 

\begin{defi}[$U(1)_c=\pm\infty$]\label{def_U1infty}
If $U(x)=\sum_{n\ge0}u_nx^n\in\R[[x]]$ and
$${\textstyle
\lim_{n\to\infty}
\sum_{j=0}^n u_j=+\infty\,,
}
$$
we say that the sum $U(1)_c=+\infty$. We similarly define that $U(1)_c=-\infty$.  
\end{defi}
So for $U(x)$ with real coefficients, we allow as values of the sums 
$U(1)_c$ also the infinities $-\infty$ and $+\infty$.

\begin{prop}\label{prop_FPS1}
Let $U(x)=\sum_{n\ge0}u_n x^n$ and $V(x)=\sum_{n\ge0}v_n x^n$ 
be in $\C[[x]]$, and let $\xi\in\C$ with $\xi\ne0$. 
If formally 
$U(x)=V(x)\frac{1}{1-\xi x}$ and the conditional sum $(V)_{1/\xi}(1)_c$ exists, then numericly
$$
\lim_{n\to\infty}\xi^{-n}u_n=(V)_{1/\xi}(1)_c\,.
$$
\end{prop}
\proof
Let $n\in\N_0$.
Since $\frac{1}{1-\xi x}=1+\xi x+\xi^2x^2+\ds$, we have 
$$
u_n=v_0\xi^n+v_1\xi^{n-1}
+\ds+v_n\xi^0=
v_0\cdot\xi^n+v_1\xi^{-1}\cdot\xi^{n}+\ds+v_n\xi^{-n}\cdot\xi^n\,.
$$
Division by $\xi^n$ and limit transition for $n\to\infty$ yield the stated formula. 
\eproof

\subsection{Three Abel's theorems}

We manipulate conditional sums $U(1)_c$ by means of three Abel's theorems, which we 
state only for the radius of convergence $\ge1$. 
The first theorem is classical, see \cite{abel,abel_Wiki} and \cite[Chapter 
II.7]{tene}.

\begin{thm}[Abel~1]\label{thm_abel}
Let
$U(x)
=\sum_{n\ge0}u_nx^n\in\C[[x]]$. 
If the conditional sum $U(1)_c$ exists, so that $U(x)\in H_1$ by  Proposition~\ref{prop_ConpowSer}, then
$$
\lim_{x\nearrow1}U(x)=U(1)_c\ \ (\in\C)\,.
$$
\end{thm}

Second Abel's theorem is seldom mentioned, but see 
\cite{abel_Wiki} for a~simple 
proof. It is given there only in the case $U(1)_c=+\infty$, but it
works for the sum and limit $-\infty$ as well.

\begin{thm}[Abel~2]\label{thm_abel2}
Let
$U(x)
=\sum_{n\ge0}u_nx^n\in
\R[[x]]\cap H_1$. If $U(1)_c=\pm\infty$ then
$$
\lim_{x\nearrow1}U(x)=U(1)_c\ \ 
(\in\{-\infty,\,+\infty\})\,.
$$
\end{thm}

Third Abel's theorem allows only nonnegative real coefficients, but it has two advantages over the two previous theorems. It applies to every formal power series  
$U(x)\in\R_{\ge0}[[x]]\cap H_1$, and it is an equivalence: one can go not 
only from the sum to the limit, but also in the reverse. We omit the 
easy proof.

\begin{thm}[Abel~3]\label{thm_abel3}
Let $U(x)
=\sum_{n\ge0}u_nx^n\in
\R_{\ge0}[[x]]\cap H_1$. Then the following limit and sum always exist, 
are possibly $+\infty$, and are equal:
$$
\lim_{x\nearrow1}U(x)=U(1)_c\ \ (\in[0,\,+\infty)\cup\{+\infty\})\,.
$$
\end{thm}
The three Abel's theorems are mutually 
incomparable. 

\subsection{Identities for conditional sums}

To reiterate our remark from the beginning of Section~\ref{sec_serGF}, we believe that some of the 
precise statements, and systematic and rigorous proofs of conversions of 
formal identities to 
numeric ones at $x=1$ in this 
and the next subsection are our original contribution.  

\begin{prop}\label{prop_totezPlus}
Let $\al,\be\in\C$ and  
$U(x), V(x), W(x)$ be in $\C[[x]]$. If formally
$$
U(x)=\al V(x)+\be W(x)
$$ 
and if the conditional sums $V(1)_c$ and $W(1)_c$ exist, then numericly
$$
U(1)_c=\al V(1)_c+\be W(1)_c\,.
$$
\end{prop}
\proof
Let $U(x)=\sum_{n\ge0}u_nx^n$, $V(x)=\sum_{n\ge0}v_nx^n$, and 
$W(x)=\sum_{n\ge0}w_nx^n$. Then
$$
{\textstyle
\sum_{j=0}^n u_j=\sum_{j=0}^n (\al v_j+\be w_j)=\al\sum_{j=0}^n v_j+\be\sum_{j=0}^n w_j\,.
}
$$
Taking the limit for $n\to\infty$, we get $U(1)_c=\al V(1)_c+\be W(1)_c$.
\eproof

\noindent
The phrase ``then numericly $U(1)_c=\al V(1)_c+\be W(1)_c$''  abbreviates the phrase ``then the conditional sum $U(1)_c$ exists and $U(1)_c=\al V(1)_c+\be 
W(1)_c$''. Sometimes we use such abbreviations. This result easily  generalizes to 
linear combinations with more than two terms. 

\begin{prop}\label{prop_Abel1}
Let $U(x)$, $V_1(x)$, $\ds$, $V_n(x)\in\C[[x]]$ be $n+1\ge2$ formal power series. If formally 
$${\textstyle
U(x)=\prod_{i=1}^n V_i(x)
}
$$ 
and if the conditional sums $U(1)_c$, $V_1(1)_c$, $\ds$, 
$V_n(1)_c$ exist, then numericly
$${\textstyle
U(1)_c=\prod_{i=1}^n V_i(1)_c\,.
}
$$
\end{prop}
\proof
Due to Propositions~\ref{prop_grouSer}, \ref{prop_prodKser}, and 
\ref{prop_ConpowSer}, and to the assumptions, we have the numeric equality  
$${\textstyle
U(x)=\prod_{i=1}^n V_i(x)\,\text{ for every number }\,x\in\C_{<1} 
}
$$
(see the proofs of Propositions~\ref{prop_bezAbel1} and 
\ref{prop_bezAbel2} for details of similar arguments). Then
\begin{eqnarray*}
U(1)_c&=&\lim_{x\nearrow1}U(x)=
\lim_{x\nearrow1}{\textstyle\prod_{i=1}^n V_i(x)}\\
&=&{\textstyle\prod_{i=1}^n\lim_{x\nearrow1}V_i(x)=\prod_{i=1}^n V_i(1)_c}\,.
\end{eqnarray*}
The first equality follows from the assumptions and 
Theorem~\ref{thm_abel}. The second equality  follows from the 
above numeric equality. The fourth and third equalities follow from the arithmetic of 
functional limits, the assumptions, and Theorem~\ref{thm_abel}. 
\eproof

\noindent
Since we have to assume the existence of $U(1)_c$, the case $n>2$ does not follow  
by iterating the case $n=2$. The example  in \cite{wikiCauchyProd} 
shows that the assumption of existence of $U(1)_c$ cannot be 
omitted.

\begin{prop}\label{prop_Abel2}
Let $U(x),V(x)\in\C[[x]]$. If formally 
$${\textstyle
U(x)=\frac{1}{V(x)} 
}
$$
and if the conditional sums $U(1)_c$ and $V(1)_c$ exist, then $V(1)_c\ne0$ and numericly
$$
U(1)_c=\frac{1}{V(1)_c}\,. 
$$
\end{prop}
\proof
Formally, $U(x)V(x)=1$. Proposition~\ref{prop_Abel1} provides 
the numeric equality $U(1)_cV(1)_c=1$.
Thus, $V(1)_c\ne0$. We divide by $V(1)_c$ and obtain the stated formula. 
\eproof

\begin{prop}\label{prop_Abel3}
Let $U(x),V(x),W(x)\in\C[[x]]$. If formally 
$${\textstyle
U(x)^2=1-\frac{V(x)}{W(x)}\,,
}
$$
if the conditional sums $U(1)_c$, $V(1)_c$, and $W(1)_c$ exist, and if $W(1)_c\ne0$, then numericly
$$
\big(U(1)_c\big)^2=1-\frac{V(1)_c}{W(1)_c}\,.
$$
\end{prop}
\proof
Formally,  
$U(x)^2\cdot W(x)=W(x)-V(x)$. 
Let $T(x)=W(x)-V(x)$. By Proposition~\ref{prop_totezPlus} we have
$T(1)_c=W(1)_c-V(1)_c$.
By Proposition~\ref{prop_Abel1}, 
$$
(U(1)_c)^2\cdot W(1)_c=T(1)_c=W(1)_c-V(1)_c\,.
$$
We divide by $W(1)_c$ and obtain the stated formula.
\eproof

\subsection{Identities for absolute sums}

We obtain analogs of the results in the previous subsection for absolute 
sums. 

\begin{prop}\label{prop_totezPlusAK}
Let $\al,\be\in\C$ and  
$U(x),V(x),W(x)\in\C[[x]]$. If formally 
$$
U(x)=\al V(x)+\be W(x)
$$ 
and if the absolute sums $V(1)_a$ and $W(1)_a$ exist, then numericly
$$
U(1)_a=\al V(1)_a+\be W(1)_a\,.
$$
\end{prop}
\proof
Let $U(x)=\sum_{n\ge0}u_nx^n$, $V(x)=\sum_{n\ge0}v_nx^n$, and 
$W(x)=\sum_{n\ge0}w_nx^n$. Since $V(1)_a$ and $W(1)_a$ are 
sums of absolutely convergent series
$\sum_{n\in\N_0}v_n$ and $\sum_{n\in\N_0}w_n$, respectively, 
and $u_n=\al v_n+\be w_n$, the stated formula follows from Proposition~\ref{prop_LKseries}.
\eproof

\begin{prop}\label{prop_bezAbel1}
Let $U(x)=\sum_{n\ge0}u_nx^n$, $V(x)=\sum_{n\ge0}v_nx^n$, and 
$W(x)=\sum_{n\ge0}w_nx^n$ be in $\C[[x]]$. If formally 
$$
U(x)=V(x)W(x)
$$ 
and if the absolute sums
$V(1)_a$ and $W(1)_a$ exist, then
the following holds.
\begin{enumerate}
\item The absolute sum $U(1)_a$ exists and numericly
$U(1)_a=V(1)_aW(1)_a$.
\item We have the bound on sums $\sum_{n\ge0}|u_n|\le\sum_{n\ge0}|v_n|\cdot\sum_{n\ge0}|w_n|$.
\end{enumerate}
\end{prop}
\proof 
1. We consider the series
$${\textstyle
S=\sum_{\langle j,\,k\rangle\in\N_0^2
}v_jw_k=\sum_{j\in\N_0}v_j\otimes
\sum_{k\in\N_0}w_k
\,.
}
$$
By Proposition~\ref{prop_proSer}, $S$ absolutely converges and has sum $V(1)_aW(1)_a$. On the other hand, $S$ has the grouping
$${\textstyle
S_P=\sum_{n\ge0}\sum_{j=0}^n v_jw_{n-j}=\sum_{n\ge0}u_n\,.
}
$$ 
Here, the inner symbol $\sum$ means sum, and the outer symbol $\sum$ means series.
By Proposition~\ref{prop_grouSer}, the series $S$ has sum $u_0+u_1+\ds=U(1)_a=V(1)_aW(1)_a$.

2. This follows from the inequality that
$$
{\textstyle
\sum_{m=0}^n|u_m|=
\sum_{m=0}^n\big|\sum_{j=0}^m v_jw_{m-j}\big|\le
\sum_{j=0}^n|v_j|\cdot\sum_{j=0}^n |w_j|
}
$$
for every $n\in\N_0$.
\eproof

\noindent
There is a~straightforward extension of this proposition to 
products with more than two factors.  

For nonnegative real coefficients, when conditional and 
absolute sums coincide, we have the following result. 

\begin{prop}\label{pridat}
Let $U(x)=\sum_{n\ge0}u_nx^n$, $V(x)=\sum_{n\ge0}v_nx^n$, and $W(x)=\sum_{n\ge0}w_nx^n$ be in $\R_{\ge0}[[x]]$, let $U(x)\ne0$, and let formally
$$
U(x)=V(x)W(x)\,. 
$$
Then the conditional sums $U(1)_c$, $V(1)_c$, and $W(1)_c$ exist, lie in $(0,+\infty)\cup\{+\infty\}$, and the numeric equality
$$
U(1)_c=V(1)_c\cdot W(1)_c\ \ (\in(0,\,+\infty)\cup\{+\infty\})
$$
holds. 
\end{prop}
\proof
The three sums exist because they are limits of non-decreasing real 
sequences. Moreover, these sums are positive. Let  $V(x)\not\in H_1$ or $W(x)\not\in 
H_1$. Then $V(1)_c=+\infty$ or $W(1)_c=+\infty$. Since 
$u_n=\sum_{i=0}^n v_iw_{n-i}$ for every $n\in\N_0$ and $v_j,w_k>0$ for 
some $j,k\in\N_0$, we deduce that $U(1)_c=u_0+u_1+\ds=+\infty$ and the
stated equality holds as 
$$
+\infty=c\cdot(+\infty)\,\text{ or }\,+\!\infty=(+\infty)\cdot c,\,\text{ with }\,c\in(0,\,+\infty)\cup
\{+\infty\}\,. 
$$
Let $V(x),W(x)\in H_1$. It follows from 
Propositions~\ref{prop_grouSer} and \ref{prop_proSer} that  $U(x)\in H_1$ 
as well, and that for every $x\in\C_{<1}$ we have the numeric identity $U(x)=V(x)W(x)$. We can therefore compute that
$$
U(1)_c=\lim_{x\nearrow1}
U(x)=\lim_{x\nearrow1}(V(x)W(x))=\lim_{x\nearrow1}V(x) 
\lim_{x\nearrow1}W(x)=
V(1)_cW(1)_c\,.
$$
The first equality 
follows from Theorem~\ref{thm_abel3}. The second equality follows from the 
above numeric identity. The fourth and third equalities follow from the arithmetic of 
functional limits, from 
Theorem~\ref{thm_abel3}, and from the fact that, due to the assumptions, 
neither of the limits is $0$ and we avoid the indefinite products $0\cdot(+\infty), 
(+\infty)\cdot0$.  
\eproof

\noindent
Without the assumption that $U(x)\ne0$, we do not always 
get the numeric equality; although formally 
$$
0+0x+0x^2+\ds=(0+0x+0x^2+\ds)(1+1x+1x^2+\ds)\,, 
$$
the numeric equality $0=0\cdot(+\infty)$ does not hold. 

Proposition~\ref{prop_Abel1} for $n=2$ and part~1 of 
Proposition~\ref{prop_bezAbel1} are incomparable. The former only 
assumes existence of the weaker conditional sums, but for all three 
generating functions. The latter assumes existence of 
the stronger absolute sums, but only for two generating functions. Yet another 
incomparable version of the product identity is 
due to F.~Mertens.

\begin{thm}[Mertens]\label{thm_Metens}
Suppose that  
$U(x)$, $V(x)$, and $W(x)$ are in $\C[[x]]$. If formally $U(x)=V(x)W(x)$, if 
the conditional sum $V(1)_c$ exists, and if the absolute sum $W(1)_a$ exists, then numericly 
$$
U(1)_c=V(1)_c\cdot W(1)_a\,.
$$    
\end{thm}
See \cite{wikiCauchyProd} for a~proof.

\begin{prop}\label{prop_bezAbel2}
Let $U(x)=\sum_{n\ge0}u_nx^n$ and $V(x)=\sum_{n\ge0}v_nx^n$ be in $\C[[x]]$, and let $v_0=1$. If formally
$${\textstyle
U(x)=\frac{1}{V(x)}
}
$$ 
and if the sum $c:=|v_1|+|v_2|+\ds<1$, 
then the following holds.
\begin{enumerate}
\item 
The absolute sums $U(1)_a$ and $V(1)_a$ exist, $V(1)_a\ne0$, and numericly $U(1)_a=
\frac{1}{V(1)_a}$.
\item The sum  
$\sum_{n\ge0}|u_n|\le\frac{1}{1-c}$.
\end{enumerate}
\end{prop}
\proof
1. $V(1)_a$ exists by the assumption on $c$. From $|V(1)_a-1|\le c<1$ we deduce that $V(1)_a\ne0$. We consider the series
$$
{\textstyle 
S=1+\sum_{\substack{m\in\N\\j_{m,1},\,j_{m,2},\,\ds,\,j_{m,m}\in\N}}(-1)^m v_{j_{m,
1}}v_{j_{m,2}}\ds v_{j_{m,m}}
}
$$
and show that it absolutely converges. Let $T$ be any finite 
subseries of $S$. Then the sum of absolute values 
of summands in $T$ is at most
$${\textstyle
1+\sum_{m\ge1}(|v_1|+|v_2|+\ds)^m=
1+\sum_{m\ge1}c^m=\frac{1}{1-c}\,.
}
$$
Thus $S$ absolutely converges.  

Let $m\in\N$. By Proposition~\ref{prop_prodKser}, the series
$${\textstyle
S_m=\bigotimes_{i=1}^m(v_1+v_2+\ds)=
\sum_{\langle
j_1,\,j_2,\,\ds,\,j_m\rangle\in\N^m}
v_{j_1}v_{j_2}\ds v_{j_m}
}
$$
absolutely converges and has 
sum $(V(1)_a-1)^m$. Thus, by Propositions~\ref{prop_LKseries} 
and \ref{prop_grouSer}, the series $S$ has the sum 
$${\textstyle
1+\sum_{m\ge1}(-1)^m(V(1)_a-1)^m=\frac{1}{V(1)_a}\,.
}
$$
On the other hand, the formal identity
$$
{\textstyle
U(x)=1+\sum_{m\ge1}(-1)^m(v_1x+v_2x^2+\ds)^m
}
$$
holds. Thus $u_0=1$ and for every $n\in\N$, 
$$
{\textstyle
u_n=\sum_{m\ge1}(-1)^m\sum_{\substack{j_{m,1},\,j_{m,2},\,\ds,\,j_{m,m}\in\N\\ j_{m,1}+j_{m,2}+\ds+j_{m,m}=n}}v_{j_{m,
1}}v_{j_{m,2}}\ds v_{j_{m,m}}\,.
}
$$
Therefore, the series $S$ has the grouping $\sum_{n\ge0}u_n$ and
by Proposition~\ref{prop_grouSer}, $S$ has the sum 
$U(1)_a=\frac{1}{V(1)_a}$.

2. This follows from the inequality that for every $n\in\N_0$, 
\begin{eqnarray*}
&&{\textstyle
\sum_{k=0}^n|u_k|=}\\
&&
{\textstyle
=1+\sum_{k=1}^n\big|\sum_{m\ge1}(-1)^m\sum_{\substack{j_{m,1},\,j_{m,2},\,\ds,\,j_{m,m}\in\N\\ j_{m,1}+j_{m,2}+\ds+j_{m,m}=k}}v_{j_{m,
1}}v_{j_{m,2}}\ds v_{j_{m,m}}\big|
}\\
&&{\textstyle
\le\sum_{m\ge0}(|v_1|+|v_2|+\ds)^m=
\frac{1}{1-c}\,.
}
\end{eqnarray*}
\eproof

\noindent
Again, Proposition~\ref{prop_Abel2} and 
item~1 of Proposition~\ref{prop_bezAbel2} are incomparable. 

The analog of Proposition~\ref{prop_Abel3} for 
absolute sums needs two lemmas.

\begin{lem}\label{lem_odmocFPS}
Let $V(x)\in\C[[x]]$ and let formally
$$
V(x)^2=tx^{2l}(1+t_1x+t_2x^2+\ds)\,,
$$
where $l\in\N_0$ and $t,t_j\in\C$;  every square in $\C[[x]]$ has 
this form. Then there exists a~number $v\in\C$ such that $v^2=t$ and that  formally  
$$
{\textstyle
V(x)=vx^l\sum_{n\ge0}\binom{1/2}{n}(t_1x+t_2x^2+\ds)^n\,.
}
$$
\end{lem}
\proof
Let $n\in\N_0$ and $a,b$ be formal variables. Vandermonde's 
identity \cite{vdmonde} states that 
$${\textstyle
\sum_{j=0}^n\binom{a}{j}\binom{b}{n-j}=\binom{a+b}{n}\ \ (\in\Q[a,\,b])\,.
}
$$
Let $V_0(x)=\sum_{n\ge0}\binom{1/2}{n}(t_1x+t_2x^2+\ds)^n$ and let 
$v_0\in\C$ be such that $(v_0)^2=t$. The instance $a=b=\frac{1}{2}$ of 
Vandermonde's identity gives, formally,
\begin{eqnarray*}
\big(v_0x^lV_0(x)\big)^2&=&{\textstyle 
tx^{2l}\sum_{n\ge0}\big(\sum_{j=0}^n\binom{1/2}{j}\binom{1/2}{n-j}\big)(t_1x+t_2x^2+\ds)^n}\\
&=&{\textstyle tx^{2l}\sum_{n\ge0}\binom{1}{n}(t_1x+t_2x^2+\ds)^n}\\
&=&
tx^{2l}(1+t_1x+t_2x^2+\ds)=V(x)^2\,.
\end{eqnarray*}
Hence $V(x)=vx^lV_0(x)$ with $v\in\{-v_0,v_0\}$.
\eproof

\noindent
We remark at this moment that (a~simple modification of) 
this argument also proves, rigorously and formally, that the unique solution $C(x)=\sum_{n\ge0}C_nx^n$ in $\C[[x]]$ of the equation $xC(x)^2-C(x)+1=0$ is
$$
{\textstyle
C(x)=\frac{1}{2x}\big(
1-\sum_{n\ge0}\binom{1/2}{n}(-4)^nx^n
\big)\,.
}
$$
This formula is invoked in some form in countless expositions on the Catalan numbers 
$C_n$, but, to our knowledge, always without a~proper rigorous formal proof (not using analytic expansions). 

\begin{lem}\label{lem_odmocU2}
Let $U(x),V(x)\in\C[[x]]$, let  formally 
$$
U(x)=V(x)^2=tx^{2l}(1+t_1x+t_2x^2+\ds)\,,
$$
where $l\in\N_0$ and $t,t_j\in\C$, and let $c:=|t_1|+|t_2|+\ds<1$. Then the absolute sums $U(1)_a$ and $V(1)_a$ exist and numericly
$$
U(1)_a=\big(V(1)_a\big)^2\,.
$$
\end{lem}
\proof
$U(1)_a$ exists by the assumption on $c$. By Lemma~\ref{lem_odmocFPS}, 
there is a~number $v\in\C$ such that 
$v^2=t$ and
$$
{\textstyle
V(x)=\sum_{n\ge0}v_nx^n=vx^l\sum_{n\ge0}\binom{1/2}{n}(t_1x+t_2x^2+\ds)^n\,.
}
$$
Thus $v_n=0$ for $n<l$, $v_l=v$, and for $n>l$, 
$$
{\textstyle
v_n=v\sum_{m\ge1}\binom{1/2}{m}\sum_{\substack{j_{m,1},\,j_{m,2},\,\ds,\,j_{m,m}\in\N\\ j_{m,1}+j_{m,2}+\ds+j_{m,m}=n-l}}t_{j_{m,
1}}t_{j_{m,2}}\ds t_{j_{m,m}}\,.
}
$$
We show that the series
$$
{\textstyle
S=0+0+\ds+0+v+\sum_{\substack{m\in\N\\j_{m,1},\,j_{m,2},\,\ds,\,j_{m,m}\in\N}}v\binom{1/2}{m} t_{j_{m,1}}t_{j_{m,2}}\ds t_{j_{m,m}}\,,
}
$$
which begins with $l$ zero summands, absolutely converges. Let $T$ be a~finite subseries of $S$. Since 
$\big|\binom{1/2}{m}\big|\le1$ for every $m\in\N$, the sum of absolute values of summands in $T$ is at most
$${\textstyle
|v|+|v|\sum_{m\ge1}(|t_1|+|t_2|+\ds)^m=
|v|+|v|\sum_{m\ge1}c^m
=\frac{|v|}{1-c}\,.}
$$
We see that $S$ absolutely converges.
The sum $V(1)_a$ exists because by Proposition~\ref{prop_grouSer}, $S$ 
has sum $v_0+v_1+\ds=V(1)_a$. 
By part~1 of Proposition~\ref{prop_bezAbel1}, $U(1)_a=\big(V(1)_a\big)^2$.
\eproof

\begin{prop}\label{prop_bezAbel3}
Let $U(x)$, $V(x)=\sum_{n\ge0}v_nx^n$, and $W(x)=\sum_{n\ge0}w_nx^n$ be in $\C[[x]]$. Suppose that the absolute sum $V(1)_a$ exists, $w_0=1$, that formally
$${\textstyle
U(x)^2=1
-\frac{V(x)}{W(x)}\,, 
}
$$
and that  
$c=|w_1|+|w_2|+\ds<1$. If $V(x)=W(x)$, then $U(x)=0$. 
Else, if $k\in\N_0$ is minimum such that $v_k\ne w_k$, $t:=w_k-v_k$ ($\ne0$), and if
$${\textstyle
\frac{1}{(1-c)|t|}\cdot\big(c|t|+\sum_{n>k}|w_n-v_n|\big)<1\,,
}
$$
then the absolute sums $U(1)_a$ and $W(1)_a$ exist, $W(1)_a\ne0$, and numericly
$$ 
\big(U(1)_a\big)^2=1-\frac{V(1)_a}{W(1)_a}\,.
$$
\end{prop}
\proof 
If $V(x)=W(x)$, then $U(x)^2=1-1=0$ and $U(x)=0$. Let 
$V(x)\ne W(x)$, and let $k$ and $t$ be as stated. We write
$$
{\textstyle
U(x)^2=\frac{W(x)-V(x)}{W(x)}=
tx^k(1+t_1x+t_2x^2+\ds)
}
$$
where $k=2l$ with $l\in\N_0$, and $t_j\in\C$. Let, formally,  
$${\textstyle
U_0(x)=\frac{1}{W(x)}=
\sum_{n\ge0}u_{0,n}x^n,\ U_1(x)=W(x)-V(x)}\,,
$$ 
and $U_2(x)=U_0(x)U_1(x)=\sum_{n\ge0}u_{2,n}x^n$.
By part~1 of Proposition~\ref{prop_bezAbel2}, the 
sums $U_0(1)_a$ and $W(1)_a$ exist, 
$W(1)_a\ne0$, and  
$U_0(1)_a=\frac{1}{W(1)_a}$. By 
Proposition~\ref{prop_totezPlusAK},  the sum $U_1(1)_a$ exists and 
$U_1(1)_a=W(1)_a-V(1)_a$. By part~1 of Proposition~\ref{prop_bezAbel1}, 
 the sum $U_2(1)_a$ exists and 
$U_2(1)_a=U_0(1)_a\cdot U_1(1)_a$.
By part~2 of Proposition~\ref{prop_bezAbel2},
$\sum_{n\ge0}|u_{0,n}|\le
\frac{1}{1-c}$. By part~2 of Proposition~\ref{prop_bezAbel1},
$$
{\textstyle
\sum_{n\ge0}|u_{2,n}|\le
\sum_{n\ge0}|u_{0,n}|
\cdot\sum_{n\ge k}|w_n-v_n|
\le\frac{1}{1-c}\big(|t|+\sum_{n>k}|w_n-v_n|\big)
\,.
}
$$
Since $U(x)^2=U_2(x)$, we have the bound
$${\textstyle
|t|(1+|t_1|+|t_2|+\ds)\le\frac{1}{1-c}\big(|t|+\sum_{n>k}|w_n-v_n|\big)\,.
}
$$
An algebraic rearrangement gives
$${\textstyle
|t_1|+|t_2|+\ds\le
\frac{1}{(1-c)|t|}\cdot\big(c|t|+\sum_{n>k}|w_n-v_n|\big)<1\,.
}
$$
By Lemma~\ref{lem_odmocU2},  the sum $U(1)_a$ exists and
$${\textstyle
\big(U(1)_a\big)^2=U_2(1)_a=
U_0(1)_a\cdot U_1(1)_a=
1-\frac{V(1)_a}{W(1)_a}\,.
}
$$
\eproof

\subsection{Generating functions of walks for convex weights}

In Subsection~\ref{subsec_PTvSMv0}, we used that $h(D_n)=1$ in every graph $G_d$ for every $n\in\N_0$. We generalize it. 

\begin{prop}\label{prop_onD}
Let $\xi\in\C$ and let $h\cc\N_2\to\C$ be a~$\xi$-convex 
light weight. Then $h(D_n)=\xi^n$ for every $n\in\N_0$, with $0^0:=1$. Thus, formally,    
$${\textstyle
\sum_{n\ge0}h(D_n) x^n=\frac{1}{1-\xi x}\ \ (\in\C[[x]])\,.
}
$$ 
In particular, for $\xi=0$ we have $\sum_{n\ge0}h(D_n) x^n=1x^0+0x^1+0x^2+\ds=1$.
\end{prop}
\proof
We prove by induction on $n\in\N_0$ that 
$h(D_n)=\xi^n$. 
For $n=0$, it holds: $h(D_0)=1$ because $D_0=\{\langle 1\rangle\}$. Let $n>0$. Then, indeed, 
\begin{eqnarray*}
h(D_n)
&=&{\textstyle\sum_{\substack{w\in D_{n-1}\\\ell(w)\in 
V(h)}}h(w)\cdot\sum_{v\in\N\setminus
\{\ell(w)\}}h(\langle\ell(w),\,v\rangle)}\\
&=&{\textstyle\sum_{w\in D_{n-1}}h(w)\cdot\xi=h(D_{n-1})\cdot\xi=\xi^{n-1}\cdot\xi=
\xi^n\,,
}
\end{eqnarray*}
where $\ell(w)$ 
denotes the last vertex of the walk $w$.
In the first equality, we take the 
walks $w'\in D_n$ such that the penultimate vertex is in $V(h)$, 
decompose them as $w'=w\,\ell(w')$, and use the 
multiplicative definition of $h(w')$. 
We can restrict to the walks $w'$: adding or dropping the condition 
$\ell(w)\in V(h)$ does not affect
the sum because walks not satisfying it have zero weight. We group
the series $\sum_{w'}h(w')$ according to the initial segment $w$, write $h(w')=
h(w)\cdot h(\langle\ell(w),\ell(w')\rangle)$, and use Propositions~\ref{prop_LKseries} and 
\ref{prop_grouSer}. Both symbols $\sum$ mean sum. In the second equality, we use the 
assumption of $\xi$-convexity and
drop the condition $\ell(w)\in V(h)$.
In the third equality, we use the other notation for the sum $\sum_{w\in 
D_{n-1}}h(w)$. In the fourth equality, we use 
induction. The last, fifth equality is trivial.  
\eproof

\section{General P\'olya's theorems for $v=1$}\label{sec_extTwoPol}

In this section, we obtain the first group of our main results, Theorems~\ref{thm_v1gen1}--\ref{thm_v1conInf}.
They generalize the formula (\ref{eq_oneMinus}) to weighted walks 
in $\vec{K}_{\N}$; the vertex to be visited coincides with 
the starting vertex. We first introduce three auxiliary 
generating functions $F_{B,h}(x)$, $F_{C,h}(x)$, and $F_{D,h}(x)$. In 
Proposition~\ref{prop_12rel} we prove two formal identities for them. 
Subsections~\ref{subsec_v=1G} and \ref{subsec_v=1C} treat general weights and $\xi$-convex weights, respectively. 

\subsection{Auxiliary generating functions}

Let $h\cc\N_2\to\C$ be a~light weight (see Definition~\ref{def_light}). 
Recall than $D_n$ is the set of walks with length $n$ in $\vec{K}_{\N}$ 
starting ar the vertex $1$.
For any length $n\in\N_0$ 
and any set $X\sus D_n$ of walks, we denote by $h(X)$ the sum of the series
$\sum_{w\in X}h(w)$. As we know from Subsection~\ref{subsec_genPT}, the main generating function is
$${\textstyle
F_{A,\,h}(x)=F_{A,\,1,\,h}(x)=
\sum_{n\ge0}h(A_n)x^n\ \ (\in\C[[x]])\,,
}
$$ 
where $A_n$ is the set of walks of length $n$ in $\vec{K}_{\N}$ that start at the vertex $1$ and later revisit it. The 
``largest'' auxiliary generating function is $F_{D,\,h}(x)=\sum_{n\ge0}h(D_n)x^n$. Another two auxiliary generating functions are
$${\textstyle
F_{B,\,h}(x)=\sum_{n\ge0}h(B_n)x^n\,\text{ and }\,
F_{C,\,h}(x)=\sum_{n\ge0}h(C_n)x^n\,,
}
$$
where $B_n$ is the set of closed walks with length $n$ in 
$\vec{K}_{\N}$, which start and end at the vertex $1$, and $C_n\sus B_n$ 
is the subset of closed walks with all inner vertices different 
from $1$. We have $h(A_0)=h(C_0)=0$ and $h(B_0)=h(D_0)=1$. For
complex nonnegative weights, we will also consider walks starting at 
another vertex $u\in\N\setminus\{1\}$. Then, for the new starting vertex $u$ 
in place of $1$, we define analogous sets 
of walks $A_n',\ds,D_n'$ and analogous generating functions 
$F_{A',h}(x),\ds,F_{D',h}(x)$; we say that they are {\em defined relative 
to the starting vertex $u$}. It follows from 
Proposition~\ref{prop_moveStVer} that the coefficients of these generating 
functions are well defined for every vertex $u\in V(h)\setminus\{1\}$.

As mentioned 
earlier, in the next proposition we generalize a~fragment of the symbolic 
method in enumerative combinatorics. It is a~new result with a~new proof. 

\begin{prop}\label{prop_12rel}
Let $h\cc\N_2\to\C$ be a~light weight. Then the above generating functions satisfy two formal identities.
\begin{enumerate}
\item $F_{A,h}(x)=F_{C,h}(x)F_{D,h}(x)$. 
\item $F_{B,h}(x)=\frac{1}{1-F_{C,h}(x)}$, equivalently, $F_{C,h}(x)=1-\frac{1}{F_{B,h}(x)}$.
\end{enumerate}
The same formal identities hold for the generating functions $F_{A',h}(x)$, $\ds$, $F_{D',h}(x)$ defined 
relative to any starting vertex $u\in V(h)\setminus\{1\}$.
\end{prop}
\proof
1. We show that for every $n\in\N_0$, the identity 
$${\textstyle
h(A_n)=\sum_{j=0}^nh(C_j)\cdot h(D_{n-j})
}
$$
holds. Let $n\le1$. Then the identity
holds because $A_n=C_n=\emptyset$ and $h(A_n)=h(C_n)=0$. For $n\ge2$, the set $A_n$ is infinite and countable. We split every walk $w\in A_n$ at the first revisit of $1$ and get a~map
$${\textstyle
\Phi\cc A_n\to
\bigcup_{j=0}^n C_j\times D_{n-j}
}
$$
defined as follows. For any walk
$$
w=\langle u_0,\,u_1,\,\ds,\,u_j,\,\ds,\,u_n
\rangle\in A_n\,,
$$ 
where 
$u_0=u_j=1$, $j>0$, and $u_i\ne1$ for $i=1,2,\ds,j-1$, we set 
$$
\Phi(w)=\langle w_1,\,w_2\rangle\ \ (\in C_j\times D_{n-j})\,,
$$
where $w_1=\langle u_0,u_1,\ds,u_j\rangle$ and 
$w_2=\langle u_j,u_{j+1},\ds,u_n\rangle$. 
It is easy to see that 
$\Phi$ is a~bijection and is
multiplicative (M): for 
every weight $h\cc\N_2\to\C$ and 
every value $\Phi(w)=\langle w_1,w_2\rangle$, we have 
$h(w)=h(w_1)h(w_2)$. 
If $X$ is a~set of pairs of walks in $\vec{K}_{\N}$, we denote by $h(X)$ 
the sum of the series
$${\textstyle
\sum_{\langle w,\,w'\rangle\in X}h(w)h(w')\,,
}
$$
if the series absolutely converges.
Let $h\cc\N_2\to\C$ be a~light weight and let $n\ge2$. Using the map $\Phi$, we get
\begin{eqnarray*}
h(A_n)&\stackrel{\mathrm{(M)}}{=}&h(\Phi[A_n])\\
&\stackrel{\text{Prop.~\ref{prop_grouSer}}}{=}&{\textstyle\sum_{j=0}^n h\big(C_j\times D_{n-j}\big)}\\
&\stackrel{\text{Prop.~\ref{prop_proSer}}}{=}&{\textstyle\sum_{j=0}^n 
h(C_j)\cdot h(D_{n-j})
}\,.
\end{eqnarray*} 
Here we use the notation $\Phi[A_n]=\{\Phi(w)\cc\;w\in A_n\}$.  

2. We show that 
$${\textstyle
F_{B,\,h}(x)=\frac{1}{1-F_{C,\,h}(x)}=1+\sum_{j\ge1}(F_{C,\,h}(x))^j
\ \ (\in\C[[x]])\,.
}
$$
We recall that $h(B_0)=1$. It suffices to show that for every $n\in\N$,
$$
{\textstyle
h(B_n)=\sum_{j\ge1}\sum_{\substack{n_1,\,n_2,\,\ds,\,n_j\in\N\\n_1+n_2+
\ds+n_j=n}}h(C_{n_1})h(C_{n_2})\ds h(C_{n_j})\,.
}
$$ 
For $n=1$ it holds because 
$h(B_1)=h(C_1)=0$ (both sets of walks are empty). For 
$n\ge2$ the set $B_n$ is infinite and countable. We split every walk in 
$B_n$ at the visits of $1$ and get a~map
$${\textstyle
\Gamma\cc B_n\to
\bigcup_{j\ge1}\bigcup_{\substack{n_1,\,n_2,\,\ds,\,n_j\in\N\\n_1+n_2+\ds+n_j=n}}C_{n_1}\times C_{n_2}\times\cdots\times 
C_{n_j}
}
$$
defined as follows. For any walk 
$$
w=\langle u_0,\,\ds,\,u_{m_1},\,\ds,\,\ds,\,u_{m_{j-1}},\,\ds,\,u_n\rangle\in B_n\,,
$$ 
where $0=m_0<m_1<\ds<m_j=n$ and $u_k=1$ iff $k\in
\{m_0,m_1,\ds,m_j\}$, we set
$$
\Gamma(w)=\langle w_1,\,w_2,\,\ds,\,w_j
\rangle\ \ (\in C_{n_1}\times C_{n_2}\times\ds\times C_{n_j} ,\,n_i=m_i-m_{i-1})\,,
$$ 
where $w_i=\langle u_{m_{i-1}},u_{m_{i-1}+1},
\ds,u_{m_i}\rangle$, $i=1,2,\ds,j$. 
It is easy to see that $\Gamma$ is a~bijection and is 
multiplicative (M), as in part~1. If $X$ is a~set of 
$k$-tuples of walks (for varying $k$), we denote by $h(X)$ the sum of the series
$${\textstyle
\sum_{\langle w_1,\,w_2,\,\ds,\,w_k)\in X}h(w_1)h(w_2)\ds h(w_k)\,,
}
$$
if the series absolutely converges.

Let $h$ be a~light weight and $n\ge2$. Using the map $\Gamma$,  we get  
\begin{eqnarray*}
h(B_n)&\stackrel{\mathrm{(M)}}{=}&h(\Gamma[B_n])\\
&\stackrel{\text{Prop.~\ref{prop_grouSer}}}{=}&{\textstyle\sum_{j\ge1}\sum_{\substack{n_1,n_2,\,\ds,\,n_j\in\N\\n_1+n_2+\ds+n_j=n}}h\big(C_{n_1}\times C_{n_2}\times\cdots\times C_{n_j}\big)}\\
&\stackrel{\text{Prop.~\ref{prop_prodKser}}}{=}&{\textstyle\sum_{j\ge1}\sum_{\substack{n_1,n_2,\,\ds,\,n_j\in\N\\n_1+n_2+\ds+n_j=n}}h(C_{n_1})h(C_{n_2})\ds h(C_{n_j})}\,.
\end{eqnarray*}

The last claim in the proposition is clear: the previous proofs really do 
not depend on which vertex is the starting one, as long as it is in $V(h)$ (so that the coefficients of the primed 
generating functions are well defined by Proposition~\ref{prop_moveStVer}).
\eproof

\subsection{General weights}\label{subsec_v=1G}

In this subsection, we drop the convexity condition for weights. Recall the conditional, absolute, and infinite sums 
$U(1)_c$, $U(1)_a$, and 
$U(1)_c=\pm\infty$ in Definitions~\ref{def_U12} and 
\ref{def_U1infty}. The generating functions $F_{A,h}(x)$, $F_{B,h}(x)$, 
$F_{C,h}(x)$, and $F_{D,h}(x)$ 
were defined in the previous subsection. In the next theorem, 
sums are conditional. The existence of conditional $A$-, $B$-, $C$-, and 
$D$-sums is assumed. 

\begin{thm}[$U(1)_c$]\label{thm_v1gen1}
Let $h\cc\N_2\to\C$ be a~light weight. If the conditional sums 
$F_{A,h}(1)_c$, $F_{B,h}(1)_c$, $F_{C,h}(1)_c$, and $F_{D,h}(1)_c$ exist, then 
$F_{B,h}(1)_c\ne0$ and
$$
F_{A,\,h}(1)_c=
\Big(1-\frac{1}{F_{B,\,h}(1)_c}\Big)F_{D,\,h}(1)_c\ \ (\in\C)\,.
$$ 
\end{thm}
\proof
Identity~1 in Proposition~\ref{prop_12rel} states that $F_{A,h}(x)=F_{C,h}(x)
F_{D,h}(x)$.
By Proposition~\ref{prop_Abel1}, 
$F_{A,h}(1)_c=F_{C,h}(1)_c
F_{D,h}(1)_c$. Identity~2 in Proposition~\ref{prop_12rel} states that $1-F_{C,h}(x)=\frac{1}{F_{B,h}(x)}$. By
Proposition~\ref{prop_Abel2},
$F_{B,h}(1)_c\ne0$ and 
$${\textstyle
(1-F_{C,h})(1)_c=\frac{1}{F_{B,h}(1)_c},\,\text{ implying }\,
F_{C,h}(1)_c=1-\frac{1}{F_{B,h}(1)_c}\,.
}
$$
The stated formula follows. 
\eproof 

We obtain an analog of the previous result for absolute sums. Now it 
suffices to assume the existence of absolute $B$- and $D$-sums.

\begin{thm}[$U(1)_a$]\label{thm_v1gen2}
Let $h\cc\N_2\to\C$ be a~light weight. If $\sum_{n\ge1}|h(B_n)|<1$ 
and if the absolute sum $F_{D,h}(1)_a$ exists, then $F_{A,h}(1)_a$ and $F_{B,h}(1)_a$ exist, $F_{B,h}(1)_a\ne0$, and
$$
F_{A,h}(1)_a=
\Big(1-\frac{1}{F_{B,\,h}(1)_a}\Big)F_{D,\,h}(1)_a\ \ (\in\C)\,.
$$ 
\end{thm}
\proof
Identity~2 in Proposition~\ref{prop_12rel} states that $1-C_h(x)=\frac{1}{B_h(x)}$. By part~1 of
Proposition~\ref{prop_bezAbel2}, the sums $(1-F_{C,h})(1)_a$ and $F_{B,h}(1)_a$ exist, $F_{B,h}(1)_a\ne0$, and  
$${\textstyle
(1-F_{C,h})(1)_a=\frac{1}{F_{B,h}(1)_a},\,\text{ implying }\,
F_{C,h}(1)_a=1-\frac{1}{F_{B,h}(1)_a}\,.
}
$$
Identity~1 in 
Proposition~\ref{prop_12rel} states that $F_{A,h}(x)=F_{C,h}(x)F_{D,h}(x)$.
By part~1 of Proposition~\ref{prop_bezAbel1}, the sum $F_{A,h}(1)_a$ exists and 
$F_{A,h}(1)_a=F_{C,h}(1)_aF_{D,h}(1)_a$. The stated formula follows.
\eproof

\noindent
Note that in the proof of Theorem~\ref{thm_v1gen2}, we have to 
start with identity~2 in Proposition~\ref{prop_12rel} to 
establish first the existence of $F_{C,h}(1)_a$. In the proof of 
Theorem~\ref{thm_v1gen1}, the order of usage of the two identities in  
Proposition~\ref{prop_12rel} is not important.

In the last theorem of this subsection, the weights are complex 
nonnegative, and the conditional $D$-sum is $+\infty$. We call a~walk $w$ 
in $\vec{K}_{\N}$ a~{\em cycle} if $|w|\ge2$ and if the first and last 
vertex of $w$ coincide. If $h\cc\N_2\to\C$ is a~weight and $X\sus\N$ is a~finite set of vertices, we say 
that $X$ is {\em cycle-covered} (with respect to $h$) if $X$ is a~subset of 
the vertex set of a~cycle $w$ with $h(w)\ne0$. 

\begin{thm}[$U(1)_c=+\infty$]
\label{thm_v1genInf}
Let $\xi\in\C_1$, let   
$h\cc\N_2\to\xi\R_{\ge0}$ 
be a~light weight, and let $(F_{D,h})_{1/\xi}(1)_c=+\infty$. The following holds.
\begin{enumerate}
\item If $\{1\}$ is cycle-covered, then also
$$
(F_{A,\,h})_{1/\xi}(1)_c=+\infty\,.
$$
\item If $u\in\N\setminus\{1\}$ is a~vertex such that $\{1,u\}$ is cycle-covered, then also 
$$
(F_{A',\,h})_{1/\xi}(1)_c=+\infty\,.
$$
Recall that $F_{A',\,h}(x)=\sum_{n\ge0}h(A_n')x^n$ with $A_n'$ being the set of walks with length $n$ in $\vec{K}_{\N}$ 
that start at $u$ and later revisit $u$.
\item If the weight $h$ is in addition symmetric, then the 
assumption of cycle-cover in part~1 may be dropped, and in part~2 it may 
be stated as the assumption that $u\in V(h)\setminus\{1\}$.
\item Without the 
assumption of cycle-cover, parts~1 and~2 generally do not hold for general weights $h$.
\end{enumerate}
\end{thm}
\proof
1. Since the starting vertex $1$ is cycle-covered, it follows from the 
subseries bound that 
$\xi^{-m}h(C_m)>0$ for some $m\in\N_0$. By identity~1 in 
Proposition~\ref{prop_12rel}, $h(A_n)=\sum_{j=0}^n h(C_j)h(D_{n-
j})$ for every $n\in\N_0$. It follows that for every $n\in\N_0$,
$$
\xi^{-m-n}h(A_{n+m})\ge\xi^{-m}h(C_m)\cdot\xi^{-n}h(D_n)\,.
$$ 
Hence 
$$
{\textstyle
\sum_{j=0}^n\xi^{-m-j}h(A_{j+m})\ge 
\xi^{-m}h(C_m)\cdot\sum_{j=0}^n \xi^{-j}h(D_j)\to+\infty
}
$$
for $n\to\infty$. We see that $(F_{A,h})_{1/\xi}(1)_c=+\infty$. 

2.  Since $\{1,u\}$  is cycle-covered, there exists 
a~$1$-$u$ walk $w_0$ and a~$u$-$1$ walk $w_1$ in $\vec{K}_{\N}$ such 
that $h(w_0)\ne0\ne h(w_1)$. For $n\in\N_0$, let $D_n'$ be the set 
of walks of length $n$ in $\vec{K}_{\N}$ starting at $u$. Due 
to the walk $w_0$ and Proposition~\ref{prop_moveStVer}, 
every series $\sum_{w\in D_n'}h(w)$, $n\in\N_0$, 
absolutely converges. The generating function 
$${\textstyle
F_{D',\,h}(x)=
\sum_{n\ge0}h(D_n')x^n
}
$$ 
is therefore well defined. We show that 
$(F_{D',h})_{1/\xi}(1)_c=+\infty$. Let $m=|w_1|$ ($\in\N_0$) and 
$c=\xi^{-m}h(w_1)>0$. For every walk $w\in D_n$, we consider the concatenated walk $w'=w_1w\in 
D_{m+n}'$. The subseries bound and Proposition~\ref{prop_LKseries} yield 
the lower bound
$$
\xi^{-m-n}h(D_{m+n}')\ge c\cdot 
\xi^{-n}h(D_n)\,.
$$
As before, since $(F_{D,h})_{1/\xi}(1)_c=+\infty$, also $(F_{D',h})_{1/\xi}(1)_c=+\infty$.
If we consider the generating function $F_{C',h}(x)$ of 
walks starting and ending at $u$ and avoiding $u$ in between, by Proposition~\ref{prop_12rel} we have the formal identity
$$
F_{A',\,h}(x)=F_{C',\,h}(x)F_{D',\,h}(x)\,.
$$
Since $(F_{D',h})_{1/\xi}(1)_c=+\infty$ and 
$\{u\}$ is cycle-covered, we deduce by part~1 that 
$$
(F_{A',\,h})_{1/\xi}(1)_c=+\infty\,.
$$

3. Let $h$ be in addition symmetric. In the situation of 
part~1, it is clear that $h(\langle1,u\rangle)\ne0$ for some 
vertex $u\in\N\setminus\{1\}$. Since 
$$
h(\langle1,u\rangle)=h(\langle u,1\rangle)\ne0\,,
$$
the $2$-cycle $\langle1,u,1\rangle$ covers $\{1\}$. In the situation of 
part~2, let $u\ne1$ be a~vertex in $V(h)$. Then there is a~walk $w_0$ in 
$\vec{K}_{\N}$ from $1$ to $u$ such that $h(w_0)\ne0$. Since $h$ is symmetric, the concatenation 
$w_0\overline{w_0}$ of $w_0$ and its reversal form a~cycle covering 
$\{1,u\}$. 

4. We consider the weight 
$$
h\cc\N_2\to1\R_{\ge0}=\R_{\ge0}
$$ 
given by 
$h(\langle m,m+1\rangle)=1$ for $m\in\N$ and by $h(\langle 
m,n\rangle)=0$ for $m,n\in\N$ with $n\ne m,m+1$. Then
$$
F_{D,\,h}(x)=1+x+x^2+\ds
$$
and $F_{D,\,h}(1)_c=+\infty$, but 
$F_{A,h}(x)=0$ and 
$F_{A,h}(1)_c=F_{A,h}(1)_a=0$. Also, 
$F_{D',h}(1)_c=+\infty$ and 
$F_{A',h}(1)_c=F_{A',h}(1)_a=0$ for every $u\in\N\setminus\{1\}$.
\eproof

\subsection{Convex weights}\label{subsec_v=1C}

In this subsection, the convexity of weights is generalized to $\xi$-convexity. The 
next theorem generalizes the formula (\ref{eq_oneMinus}) to complex 
nonnegative $\xi$-convex weights. In fact, the proof generalizes the 
argument presented in Subsection~\ref{subsec_PTvSMv0}. We allow starting vertices $u\ne1$.
A~weight $h\cc\N_2\to\C$ is {\em probabilistic} if $h\cc\N_2\to\R_{\ge0}$ and $h$ is $1$-convex.

\begin{thm}[$U(1)_c$ and  $\xi^{-1}h\ge0$]\label{thm_v1conNN}
Let $\xi\in\C_1$ and let  $h\cc\N_2\to\xi\R_{\ge0}$ be a~$\xi$-convex light weight. The following holds.
\begin{enumerate}
\item The 
conditional sum $(F_{B,h})_{1/\xi}(1)_c$ exists with $(F_{B,h})_{1/\xi}(1)_c\in[1,+\infty)\cup\{+\infty\}$ 
and 
$$
\lim_{n\to\infty}\xi^{-n}h(A_n)=
1-\frac{1}{(F_{B,\,h})_{1/\xi}(1)_c}\in[0,\,1]\,.
$$

\item If $\{1\}$ is cycle-covered, which is the case when $h$
is symmetric, then the last interval may be replaced with $(0,1]$.

\item Let $u\in V(h)\setminus\{1\}$. Then the 
conditional sum $(F_{B',h})_{1/\xi}(1)_c$ exists with $(F_{B',h})_{1/\xi}(1)_c\in[1,+\infty)\cup\{+\infty\}$ 
and 
$$
\lim_{n\to\infty}\xi^{-n}h(A_n')=
1-\frac{1}{(F_{B',\,h})_{1/\xi}(1)_c}\in[0,\,1]\,.
$$
Here $A_n'$ and $F_{B',h}(x)$ are defined relative to the starting vertex $u$. If $\{u\}$ 
is cycle-covered, which is the case when the weight $h$
is symmetric, then the last interval may be replaced with $(0,1]$.

\item {\em (classical criterion of transience)} If $h$ is a~probabilistic weight, then 
$$
\text{$h$ is transient, i.e. }
\lim_{n\to\infty}h(A_n)<1\iff F_{B,\,h}(1)_c<+\infty\,.
$$
\end{enumerate}
\end{thm}
\proof
1. The sum $(F_{B,h})_{1/\xi}(1)_c$ exists or equals $+\infty$ due to 
non-negativity: 
$\xi^{-1}h\ge0$ implies $\xi^{-n}h(B_n)\ge0$ for every $n\in\N_0$. 
In the former case we have $(F_{B,h})_{1/\xi}(1)_c\ge1$, because 
$h(B_0)=1$. We prove the displayed identity.
Identity~2 in Proposition~\ref{prop_12rel} states that $F_{C,h}(x)=
1-\frac{1}{F_{B,h}(x)}$. Equivalently, 
$$
F_{C,\,h}(x)F_{B,\,h}(x)=F_{B,\,h}(x)-1\,.
$$
By Proposition~\ref{prop_onD}, $h(D_n)=\xi^n$ for every $n\in\N_0$.
Since $B_n,C_n\sus D_n$ for every $n\in\N_0$, the subseries bound
implies
$$
\xi^{-n}h(B_n),\,
\xi^{-n}h(C_n)\in[0,\,1]\,.
$$
Thus, the power series $(F_{B,h})_{1/\xi}(x)$ and 
$(F_{C,h})_{1/\xi}(x)$ absolutely converge for every number $x\in[0,1)$.
By Propositions~\ref{prop_LKseries}, 
\ref{prop_grouSer}, and \ref{prop_proSer}, the numeric identity
$$
(F_{C,\,h})_{1/\xi}(x)(F_{B,\,h})_{1/\xi}(x)=(F_{B,\,h})_{1/\xi}(x)-1
$$
holds for every number $x\in[0,1)$. For every such $x$, we have $(F_{B,h})_{1/\xi}(x)\ge1$. We divide by  
$(F_{B,h})_{1/\xi}(x)$ and get that
$$
{\textstyle
(F_{C,\,h})_{1/\xi}(x)=
1-\frac{1}{(F_{B,\,h})_{1/\xi}(x)}
\,\text{ for every }\,x\in[0,\,1)\,.
}
$$
Using Theorem~\ref{thm_abel3} and the arithmetic of functional limits, we 
deduce from this that
\begin{eqnarray*}
(F_{C,\,h})_{1/\xi}(1)_c&=&{\textstyle
\lim_{x\nearrow1}(F_{C,\,h})_{1/\xi}(x)=
\lim_{x\nearrow1}\big(1-\frac{1}{(F_{B,\,h})_{1/\xi}(x)}\big)}\\
&=&{\textstyle 1-\frac{1}{\lim_{x\nearrow1}(F_{B,\,h})_{1/\xi}(x)}=1-\frac{1}{(F_{B,\,h})_{1/\xi}(1)_c}
\ \ (\in[0,\,1])\,. 
}   
\end{eqnarray*}
Identity~1 in 
Proposition~\ref{prop_12rel} states that $F_{A,h}(x)=F_{C,h}(x)F_{D,h}(x)=
F_{C,h}(x)\frac{1}{1-\xi x}$.
By Proposition~\ref{prop_FPS1}, 
$$
{\textstyle
\lim_{n\to\infty}\xi^{-n}h(A_n)=(F_{C,\,h})_{1/\xi}(1)_c=1-\frac{1}{(F_{B,\,h})_{1/\xi}(1)_c}\in[0,\,1]\,.
}
$$

2. If $\{1\}$ is cycle-covered, then the subseries bound implies that 
$\xi^{-m}h(B_m)>0$ for some integer $m\ge2$, and $(F_{B,h})_{1/\xi}(1)_c>1$. If $h$ is symmetric, then 
$\{1\}$ is cycle-covered by a~$2$-cycle $\langle1,u,1\rangle$ because $h$ is $\xi$-convex with $\xi\ne0$.

3. The assumption on $u$ and Proposition~\ref{prop_moveStVer} 
ensure that the generating function
$F_{D',h}(x)$ is correctly defined.
The proof of part~1 is based on the formal identities in 
Proposition~\ref{prop_12rel}, which hold regardless of the starting 
vertex. We therefore get the analog of part~1 without any change in the proof. In part~2 we assume that 
$\{u\}$ is cycle-covered. If $h$ is symmetric, $\{u\}$ is again 
cycle-covered by a~$2$-cycle: $u\in V(h)$ and $h$ is $\xi$-convex with 
$\xi\ne0$.

4. This is an immediate corollary of part~$1$.
\eproof

\noindent
The simple example given in the proof of part~4 of 
Theorem~\ref{thm_v1genInf} shows that in part~1 of the previous theorem, we 
may have
$$
\lim_{n\to\infty}\xi^{-n}h(A_n)=0\,.
$$
The criterion of transience in part~4 of the previous theorem is well 
known. It is stated, for example, in probabilistic language in Billingsley \cite[Theorem~8.2 (i), p.~118]{bill}. To complete the proof that part~4 is really the classical 
transience criterion, it remains to show that 
$\lim_{n\to\infty}h(A_n)<1$ if and only if the random walk on $\vec{K}_{\N}$ determined by 
$h$ is transient in the (standard) sense of Markov chains, and that the 
analogous result holds for the right-hand side of the equivalence. We give
these proofs in \cite{klaz2}.

The next two theorems treat conditional and absolute sums, and 
they generalize the case $d\ge3$ of 
the formula (\ref{eq_oneMinus}) to $\xi$-convex weights. The case $d\le2$ is 
generalized in 
Theorem~\ref{thm_v1conInf}. 
\begin{quote}
{\em The next three theorems therefore generalize the classical criterion of transience in part~4 of Theorem~\ref{thm_v1conNN} to complex weights.}    
\end{quote}
The case  of $0$-convex weights is always treated 
separately. In the next theorem, we assume the existence of conditional $B$- 
and $C$-sums.  

\begin{thm}[$U(1)_c$]\label{thm_v1con1}
Let $h\cc\N_2\to\C$ be a~$\xi$-convex light weight. The following holds.
\begin{enumerate}
\item If $\xi\in\C_1$ and the conditional sums 
$(F_{B,h})_{1/\xi}(1)_c$ and $(F_{C,h})_{1/\xi}(1)_c$ exist, then   $(F_{B,h})_{1/\xi}(1)_c\ne0$ and
$$
\lim_{n\to\infty}\xi^{-n}h(A_n)=
1-\frac{1}{(F_{B,\,h})_{1/\xi}(1)_c}\ \ (\in\C)\,.
$$
\item If $\xi=0$ and the conditional sums $F_{B,h}(1)_c$ and $F_{C,h}(1)_c$ exist, then 
$F_{B,h}(1)_c\ne0$ and
$$
F_{A,\,h}(1)_c=1-\frac{1}{F_{B,\,h}(1)_c}\ \ (\in\C)\,.
$$
\end{enumerate}
\end{thm}
\proof
1. By identity~1 in 
Proposition~\ref{prop_12rel} and by  
Proposition~\ref{prop_onD},  
$F_{A,h}(x)=F_{C,h}(x)F_{D,h}(x)=F_{C,h}(x)\frac{1}{1-\xi x}$. By Proposition~\ref{prop_FPS1},
$$
\lim_{n\to\infty}\xi^{-n}h(A_n)=
(F_{C,\,h})_{1/\xi}(1)_c\,.
$$
Identity~2 of Proposition~\ref{prop_12rel} states that 
$1-F_{C,h}(x)=\frac{1}{F_{B,h}(x)}$. By Proposition~\ref{prop_Abel2} we have 
$(F_{B,h})_{1/\xi}(1)_c\ne0$ and 
$${\textstyle
(1-F_{C,\,h})_{1/\xi}(1)_c=
\frac{1}{(F_{B,\,h})_{1/\xi}(1)_c},\,\text{ so that}\,
(F_{C,\,h})_{1/\xi}(1)_c=
1-\frac{1}{(F_{B,\,h})_{1/\xi}(1)_c}\,.
}
$$
The stated formula follows.

2. Using identity~1 in 
Proposition~\ref{prop_12rel} and Proposition~\ref{prop_onD}, we get 
that $F_{A,h}(x)=C_h(x)$. We use the displayed formula for
$(F_{C,h})_{1/\xi}(1)_c$ with $\xi=1$.
\eproof

In the case of absolute sums, it suffices to assume only the existence of $B$-sums. 

\begin{thm}[$U(1)_a$]\label{thm_v1con2}
Suppose that $h\cc\N_2\to\C$ is a~$\xi$-convex light weight and that  $\sum_{n\ge1}\big|h(B_n)\big|<1$. The following holds.
\begin{enumerate}
\item If $\xi\in\C_1$, then the absolute sum $(F_{B,h})_{1/\xi}(1)_a$ exists, $(F_{B,h})_{1/\xi}(1)_a\ne0$, and
$$
\lim_{n\to\infty} \xi^{-n}h(A_n)=
1-\frac{1}{(F_{B,\,h})_{1/\xi}(1)_a}\ \ (\in\C)\,.
$$
\item If $\xi=0$, then
the absolute sums $F_{A,h}(1)_a$ and 
$F_{B,h}(1)_a$ exist, 
$F_{B,h}(1)_a\ne0$, and
$$
F_{A,\,h}(1)_a=
1-\frac{1}{F_{B,\,h}(1)_a}\ \ (\in\C)\,.
$$
\end{enumerate}
\end{thm}
\proof 
1. Identity~2 in Proposition~\ref{prop_12rel} states that 
$1-F_{C,h}(x)=\frac{1}{F_{B,h}(x)}$. By part~1 of 
Proposition~\ref{prop_bezAbel2}, 
the sums $(F_{B,h})_{1/\xi}(1)_a$ and $(1-F_{C,h})_{1/\xi}(1)_a$ exist, 
$(F_{B,h})_{1/\xi}(1)_a\ne0$, and
$${\textstyle
(1-F_{C,\,h})_{1/\xi}(1)_a=
\frac{1}{(F_{B,\,h})_{1/\xi}(1)_a}\,\text{ hence }\,
(F_{C,\,h})_{1/\xi}(1)_a=1-\frac{1}{(F_{B,\,h})_{1/\xi}(1)_a}\,.
}
$$
By identity~1 in 
Proposition~\ref{prop_12rel} and 
Proposition~\ref{prop_onD}, we have
$F_{A,h}(x)=F_{C,h}(x)F_{D,h}(x)=F_{C,h}(x)\frac{1}{1-\xi x}$.
Then
$${\textstyle
\lim_{n\to\infty}\xi^{-n}h(A_n)=(F_{C,\,h})_{1/\xi}(1)_c=(F_{C,\,h})_{1/\xi}(1)_a
=1-\frac{1}{(F_{B,\,h})_{1/\xi}(1)_a}\,.
}
$$
The first equality follows from Proposition~\ref{prop_FPS1}. The 
second equality follows from the fact that the sum $(F_{C,h})_{1/\xi}(1)_a$ 
exists. The third equality was proven above. 

2. By identity~1 in 
Proposition~\ref{prop_12rel} and by Proposition~\ref{prop_onD},   
$F_{A,h}(x)=F_{C,h}(x)$ as formal power series. We use the displayed formula for
$(F_{C,h})_{1/\xi}(1)_a$ with $\xi=1$.
\eproof

The last theorem of this subsection treats the case $U(1)_c=\pm\infty$.
We allow starting vertices $u\ne1$.

\begin{thm}[$U(1)_c=\pm\infty$]
\label{thm_v1conInf} 
Let $\xi\in\C_1$ and let $h\cc\N_2\to\xi\R$ be a~light weight. The following holds.
\begin{enumerate}
\item If $h\cc\N_2\to\xi\R_{\ge0}$, if $h$ is $\xi$-convex, and if $(F_{B,h})_{1/\xi}(1)_c=+\infty$, then
$$
\lim_{n\to\infty}\xi^{-n}h(A_n)=1\,.
$$

\item Under the assumptions of part~$1$ on $h$, if $u\ne1$ is a~vertex such that $\{1,u\}$ is cycle-covered, then also
$$
\lim_{n\to\infty}\xi^{-n}h(A_n')=1\,.
$$
Recall that $A_n'$ is the set of walks with length $n$ in $\vec{K}_{\N}$ that 
start at $u$ and later revisit $u$.
If $h$ is in addition symmetric, it suffices to assume that $u\in V(h)\setminus\{1\}$.

\item If $h$ is $0$-convex, if $F_{B,h}(x),F_{C,h}(x)\in H_1$, and if $(F_{B,h})_{1/\xi}(1)_c=
\pm\infty$, then 
$$
\lim_{x\nearrow1}(F_{A,\,h})_{1/\xi}(x)=1\,.
$$ 
\end{enumerate}
\end{thm}
\proof
1. This is case $(F_{B,h})_{1/\xi}(1)_c=+\infty$ of part~1 of Theorem~\ref{thm_v1conNN}.

2. Since $\{1,u\}$ is cycle-covered,  there exists 
a~$1$-$u$ walk $w_0$ and a~$u$-$1$ walk $w_1$ such that $h(w_0)\ne0\ne 
h(w_1)$. Thus $u\in V(h)$ and the generating function $F_{B',h}(x)$ is 
well defined. We show that 
$(F_{B',h})_{1/\xi}(1)_c=+\infty$. This follows from the assumption that 
$(F_{B,h})_{1/\xi}(1)_c=+\infty$ and from the inequality that
$$
\xi^{-m_1-m_0-n}h(B_{m_1+m_0+n}')\ge
c\cdot\xi^{-n}h(B_n)
$$
for every $n\in\N_0$, where $c=\xi^{-m_1-m_0}h(w_1)h(w_0)>0$ and 
$m_i=|w_i|$. We obtain this inequality by considering, for every 
cycle $w$ covering $\{1\}$, the cycle $w'=w_1ww_0$ that covers $\{u\}$. The 
rest follows from the case $(F_{B',h})_{1/\xi}(1)_c=+\infty$ of part~3 of Theorem~\ref{thm_v1conNN}. 
If $h$ is symmetric and $u\in V(h)\setminus\{1\}$, then 
$\{1,u\}$ is cycle-covered by $w_0\overline{w_0}$.   

3. We claim that for some $\de>0$ we have numericly 
$$
{\textstyle
\text{$(F_{C,\,h})_{1/\xi}(x)=
1-\frac{1}{(F_{B,\,h})_{1/\xi}(x)}$ for every $x\in(1-\de,\,1)$}\,.
}
$$
This follows from the assumption that 
$F_{B,h}(x),F_{C,h}(x)\in H_1$, from identity~2 in 
Proposition~\ref{prop_12rel}, from Propositions~\ref{prop_LKseries}, \ref{prop_grouSer}, and
\ref{prop_proSer}, and 
from the fact that $(F_{B,h})_{1/\xi}(x)\ne0$ for every number 
$x<1$ close to $1$. The last fact follows from the
limit
$$
\lim_{x\nearrow1}(F_{B,\,h})_{1/\xi}(x)=
(F_{B,\,h})_{1/\xi}(1)_c=\pm\infty\,,
$$
which follows from the assumption and Theorem~\ref{thm_abel2}. We see that
\begin{eqnarray*}
\lim_{x\nearrow1}
(F_{A,h})_{1/\xi}(x)&=&
\lim_{x\nearrow1}(F_{C,\,h})_{1/\xi}(x)=\lim_{x\nearrow1}
{\textstyle
\big(1-\frac{1}{(F_{B,\,h})_{1/\xi}(x)}\big)}\\
&=&{\textstyle
1-\frac{1}{\lim_{x\nearrow1}(F_{B,\,h})_{1/\xi}(x)}=
1-\frac{1}{(F_{B,\,h})_{1/\xi}(1)_c}=
1-\frac{1}{\pm\infty}=1\,.
}
\end{eqnarray*}
The first equality follows from identity~1 in Proposition~\ref{prop_12rel} and Proposition~\ref{prop_onD}. The second equality follows 
from the above claim. The fourth and third equalities follow from Theorem~\ref{thm_abel2} and
from the arithmetic of functional limits. The fifth equality follows from the assumption.
The last, sixth equality follows from the arithmetic of infinities. 
\eproof

\noindent
We reiterate that the three previous theorems are generalizations of the 
classical criterion of transience of random walks on $\N$, which is stated, for example, in  \cite[Theorem~8.2 (i), p.~118]{bill}.

\section{General P\'olya's theorems for $v\ne1$}\label{sec_extTwoPol2}

In this section, we obtain the second group of our main results, which are Theorems~\ref{thm_vne1gen1}--\ref{thm_vne1conInf}.
They generalize formulas (\ref{eq_exactly1}) and (\ref{eq_odmoc}) to the case $v\ne1$ in the edge-weighted graph $\vec{K}_{\N}$. 
We always assume that $v\in V(h)$, which means 
that there is a~$1$-$v$ walk $w$ with $h(w)\ne0$, for else the main 
generating function $F_{A,v,h}(x)=0$.
We first introduce three additional auxiliary generating functions $F_{B,v,h}(x)$, 
$F_{C,v,h}(x)$, and $F_{E,v,h}(x)$, and in 
Proposition~\ref{prop_12relvne1} we sketch proofs of 
three formal identities for them. 
In Subsections~\ref{subsec_vne1} and \ref{subsec_vne1C}, we treat general weights and $\xi$-convex weights, 
respectively. The second identity in Proposition~\ref{prop_12relvne1} 
contains a square, and in Corollary~\ref{cor_monoton} we deduce from this a~result on monotonicity. The 
transitivity of the grid graph $G_d$ plays an important role in the case $\overline{v}\ne\overline{0}$, and this is reflected 
in Definition~\ref{def_vTrans} of 
$v$-transitive weights. 
In Proposition~\ref{prop_onFv}, we present, for the sake of interest, two formal identities for 
the generating function $F_{F,v,h}(x)$ of the $1$-$v$ walks in $\vec{K}_{\N}$.

\subsection{Auxiliary generating functions}

Let $v\in\N\setminus\{1\}$ be the vertex to be visited and $h\cc\N_2\to\C$ be 
a~light weight (recall Definition~\ref{def_light}). As we know from Section~\ref{subsec_genPT}, the main generating function is
$${\textstyle
F_{A,\,v,\,h}(x)=
\sum_{n\ge0}h(A_{v,\,n})x^n\ \ (\in\C[[x]])\,,
}
$$ 
where $A_{v,n}$ is the set of walks of length $n$ 
in $\vec{K}_{\N}$ starting at $1$ and 
visiting $v$, and $h(A_{v,n})$ is the sum of the series
$\sum_{w\in A_{v,n}}h(w)$. Let $n\in\N_0$. Besides the sets of walks 
$A_{v,n}$ and the earlier defined $B_n$, $C_n$, and $D_n$, we introduce three new  
families of sets of walks: 
$B_{v,n}$ are the walks with length $n$ in $\vec{K}_{\N}$ 
that start and end at $1$ and avoid $v$, $C_{v,n}$ are the walks of length $n$ in $\vec{K}_{\N}$ that 
start at $1$, end at $v$, and have inner vertices different
from $v$, and $E_{v,n}$ are the walks with length $n$ in $\vec{K}_{\N}$ that start at $1$, end at $v$, and have inner 
vertices different from both $1$ and $v$. The corresponding auxiliary generating functions are
\begin{eqnarray*}
F_{B,\,v,\,h}(x)&=&{\textstyle\sum_{n\ge0}h(B_{v,\,n})x^n,\ F_{C,\,v,\,h}(x)=\sum_{n\ge0}h(C_{v,\,n})x^n,\,\text{ and}}\\ 
F_{E,\,v,\,h}(x)&=&{\textstyle\sum_{n\ge0}h(E_{v,\,n})x^n\,.}
\end{eqnarray*}
We have $h(A_{v,0})=h(C_{v,0})=
h(E_{v,0})=0$ and 
$h(B_0)=h(B_{v,0})=h(D_0)=1$. 
As in the case $v=1$, for complex nonnegative weights we investigate 
the effects of replacing the starting vertex $1$ with another vertex $u$, 
and we define analogous primed sets of walks and analogous primed 
generating functions. Again, by Proposition~\ref{prop_moveStVer}, the 
primed generating functions have well defined coefficients for every vertex 
$u\in V(h)\setminus\{1\}$.

We can obtain formal identities 
for the generating functions with $v\ne1$ by the extended symbolic method only
if we restrict to $v$-transitive weights, which we introduce next.

\begin{defi}[$v$-transitivity]\label{def_vTrans}
Let $v\in\N$ be a~vertex and let $h\cc\N_2\to\C$ be an edge weight. We say that $h$ is $v$-transitive if there exists 
a~bijection $f\cc\N\to\N$ such that $f(1)=v$ and 
for every edge $\langle u,u'\rangle\in\N_2$, the equality
$$
h(\langle u,\,u'\rangle)=h(\langle f(u),\,f(u')\rangle)
$$
holds.
\end{defi}
We only sketch the 
proof of the next proposition, which is analogous to Proposition~\ref{prop_12rel}. We use 
the operation of {\em split/join} of two (or more) walks. 
If $w=\langle u_0,\ds,u_n\rangle$ with $n\in\N_0$ is a~walk in 
$\vec{K}_{\N}$ and 
$i\in\{0,\ds,n\}$, we write $w=w_1\vee w_2$ for the split of $w$ 
at $u_i$ in the two walks
$$
w_1=\langle u_0,\,\ds,\,u_i\rangle\,\text{ and }\,
w_2=\langle
u_i,\,\ds,\,u_n\rangle\,.
$$
In the reverse, if $w_1$ and $w_2$ are walks in $\vec{K}_{\N}$ such that the last vertex of $w_1$ coincides with the first vertex of $w_2$, we write
$w=w_1\vee w_2$ for the join of the two walks when $w_1$ is simply 
continued by $w_2$. It is straightforward to extend these two 
binary operations to more than two 
arguments.

\begin{prop}\label{prop_12relvne1}
Let $v\in\N\setminus\{1\}$,  $h\cc\N_2\to\C$ be 
a~$v$-transitive light weight, and
$F_{A,v,h}(x)$, $F_{B,h}(x)$, $F_{B,v,h}(x)$, $F_{C,v,h}(x)$, $F_{D,h}(x)$, and $F_{E,v,h}(x)$ be the above  generating functions. They satisfy three formal identities. 
\begin{enumerate}
\item $F_{A,v,h}(x)=
F_{C,v,h}(x)F_{D,h}(x)$.
\item $F_{B,h}(x)=
F_{B,v,h}(x)+
F_{C,v,h}(x)^2F_{B,h}(x)$,
equivalently, 
$${\textstyle
F_{C,v,h}(x)^2=1-
\frac{F_{B,v,h}(x)}{F_{B,h}(x)}\,.
}
$$
\item $F_{C,v,h}(x)=
F_{B,v,h}(x)F_{E,v,h}(x)$.
\end{enumerate}
The same formal identities hold for the generating functions $F_{A',v,h}(x)$, $\ds$, $F_{E',v,h}(x)$ defined 
relative to any starting vertex $u\in V(h)\setminus\{1,v\}$, if the weight 
$h$ is in addition $u$-transitive (for identity~3 this assumption is not needed).
\end{prop}
\proof (Sketches.)
Let $v$ and $h$ be as stated. We consider the following sets of walks in $\vec{K}_{\N}$ starting at $1$: 
$$
{\textstyle
A_v=\bigcup_{n\ge0}A_{v,\,n},\ B=\bigcup_{n\ge0}B_n,\ \ds,\ E_v=\bigcup_{n\ge0}E_{v,\,n}\,.
}
$$
Let $f$ be a~bijection in Definition~\ref{def_vTrans}.

1. We split every walk $w\in A_v$ at the first visit of $v$
as $w=w_1\vee w_2$. Then 
$|w|=|w_1|+|w_2|$, $w_1\in C_v$, $f^{-1}[w_2]\in D$, and $h(f^{-1}[w_2])=h(w_2)$. The stated identity follows.

2. It suffices to prove the first identity. We consider the disjoint union
$$
B=B_v\cup B_v^*\,,
$$
where $B_v^*$ is the set of walks in $\vec{K}_{\N}$ that start and end at $1$ 
and visit the vertex $v$. We split every walk $w\in B_v^*$ at the first 
and last visits of $v$ as $w=w_1\vee w_2\vee w_3$. Then 
$|w|=|w_1|+|w_2|+|w_3|$, the walks $w_1,\overline{w_3}
\in C_v$, the walk $f^{-1}[w_2]\in B$, and $h(f^{-1}[w_2])=h(w_2)$, where $\overline{w_3}$ is the reversal of $w_3$. The stated identity follows.

3. We split every walk $w\in C_v$ at the last visit 
of $1$ as $w=w_1\vee w_2$. 
Then $|w|=|w_1|+|w_2|$, $w_1\in B_v$, and $w_2\in E_v$.
Without using the $v$-transitivity of $h$, we obtain the stated identity.

Suppose that $u\in V(h)\setminus
\{1,v\}$ is such that $h$ is $u$-transitive, which means that there 
exists a~weight-preserving bijection $f_0\cc\N\to\N$ with $f(1)=u$. 
To obtain the identities 1 and 2 for the primed generating functions  
defined relative to the starting vertex $u$, we need a~weight-preserving bijection $f_1\cc\N\to\N$ 
such that $f(u)=v$. We set 
$f_1=f(f_0^{-1})$. Identity~3 holds with the primed generating functions 
for any $u\in V(h)\setminus\{1,v\}$.
\eproof

For the sake of interest, we mention and sketch the proofs of two identities for the generating 
function   
$${\textstyle
F_{F,\,v,\,h}(x)
=\sum_{n\ge0}h(F_{v,\,n})x^n\,,
}
$$
where $F_{v,n}$ is the set of the walks of length $n$ in $\vec{K}_{\N}$ that start at $1$ and end at the vertex $v$.

\begin{prop}\label{prop_onFv}
Let $v\in\N\setminus\{1\}$ and  $h\cc\N_2\to\C$ be 
a~$v$-transitive light weight. Then two formal identities hold. 
\begin{enumerate}
\item  $F_{F,v,h}(x)=F_{C,v,h}(x)F_{B,h}(x)$.
\item $F_{B,v,h}(x)=F_{B,h}(x)-\frac{F_{F,v,h}(x)^2}{F_{B,h}(x)}$.
\end{enumerate}
\end{prop}
\proof
(Sketches.) Let $v$ and $h$ be as stated. We define the set of walks in $\vec{K}_{\N}$
$$
{\textstyle
F_v=\bigcup_{n\ge0}F_{v,\,n}\,.
}
$$
Let $f$ be a~bijection in Definition~\ref{def_vTrans}.

1. We consider the disjoint union
$$
F_v=C_v\cup F_v^*\,,
$$
where $F_v^*$ are the walks in $\vec{K}_{\N}$ that start at $1$, 
end at $v$, and visit $v$ at least twice. We split every walk $w\in F_v^*$ at the first visit of $v$
as $w=w_1\vee w_2$. Then $|w|=|w_1|+|w_2|$, 
$w_1\in C_v$, $f^{-1}[w_2]\in(B\setminus
\{\langle1\rangle\})$, and 
$h(f^{-1}[w_2])=h(w_2)$. It follows that
$$
F_{v,\,h}(x)=F_{C,\,v,\,h}(x)+F_{C,\,v,\,h}(x)(F_{B,\,h}(x)-1)=F_{C,\,v,\,h}(x)F_{B,\,h}(x)\,.
$$

2. We prove the algebraic rearrangement
$$
F_{B,\,v,\,h}(x)\cdot 
F_{B,\,h}(x)+F_{F,\,v,\,h}(x)^2=F_{B,\,h}(x)^2
$$
of the stated identity. We consider the disjoint union
$$
(B_v\times B)\cup G=B\times B\,,
$$
where $G\sus B\times B$ are the pairs of walks $\langle 
w,w'\rangle\in B\times B$ such that $w$ visits the vertex $v$. 
Clearly, the 
former Cartesian product is counted  by $F_{B,v,h}(x)\cdot F_{B,h}(x)$, and the 
latter by $F_{B,h}(x)^2$. It remains to show that 
$G$ is counted by $F_{v,h}(x)^2$. We prove it by defining two mutually inverse maps
$$
g\cc G\to F_v\times F_v\,\text{ and }\,g'\cc F_v\times F_v\to G
$$
such that if $g(\langle w_1,w_2\rangle)=\langle 
w_1',w_2'\rangle$, then $h(w_1)h(w_2)=h(w_1')h(w_2')$, and 
similarly for $g'$. As before,
we denote the reversal of a~walk by bar. Let $\langle w_1,w_2\rangle\in G$. In order to define $g$, we 
split $w_1$ at the last occurrence of $v$ as $w_1=w_1'\vee w_1''$.  
We set
$$
g(w_1,\,w_2)=\langle w_1',\,\overline{w_1''}\vee f[w_2]
\rangle\,.
$$
Let $\langle w_3,w_4\rangle\in F_v\times F_v$. In order to define $g'$, we split $w_4$ at the first occurrence of $v$ as $w_4=w_4'\vee w_4''$. We set
$$
g'(w_3,\,w_4)=\langle w_3\vee\overline{w_4'},\,f^{-1}[w_4'']
\rangle\,.
$$
It is not hard to check that the maps $g$ and $g'$ have the stated 
properties.
\eproof

Identity~2 in Proposition~\ref{prop_12relvne1}
and identity~1 in Proposition~\ref{prop_onFv} 
provide two expressions for 
$F_{C,v,h}(x)$:
$${\textstyle
F_{C,\,v,\,h}(x)^2=
1-\frac{F_{B,\,v,\,h}(x)}{F_{B,\,h}(x)}\,\text{ and }\,
F_{C,\,v,\,h}(x)=\frac{F_{F,\,v,\,h}(x)}{F_{B,\,h}(x)}\,.
}
$$
The latter expression is simpler than the former, both in 
appearance and factually (the set $F_v$ has a~simpler definition than 
the set $B_v$),
but we use the former expression. The reason is that since
$B_{v,n}\sus B_n$, the former expression implies at once in 
Theorem~\ref{thm_vne1conNN} that the root lies in $(0,1)$. The latter 
expression involves disjoint sets of  walks $F_{v,n}$ and $B_n$, which is inconvenient.

In the next corollary, we utilize the square in identity~2 in 
Proposition~\ref{prop_12relvne1}. As we know, the sum of a~subseries is 
generally unrelated to the sum of the series. However, in the case of 
inclusions $B_{v,n}\sus B_n$, $n\in\N_0$, and real 
weights $h$, there is a~corresponding inequality between the values of $F_{B,v,h}(x)$ 
and $F_{B,h}(x)$.

\begin{cor}\label{cor_monoton}
Let $v\in\N\setminus\{1\}$, $h\cc\N_2\to\R$ be a~$v$-transitive light weight, and let $F_{B,v,h}(x),F_{B,h}(x),F_{C,v,h}(x)\in H_1$. Then for every number $x\in(-1,1)$,  the following holds.
\begin{enumerate}
\item If $F_{B,h}(x)=0$ then also 
$F_{B,v,h}(x)=0$.
\item If $F_{B,h}(x)\ne 0$ then
$$
\frac{F_{B,\,v,\,h}(x)}{F_{B,\,h}(x)}\le1\,.
$$
\end{enumerate}
\end{cor}
\proof
By identity~2 in Proposition~\ref{prop_12relvne1}, 
the assumptions, and Propositions~\ref{prop_LKseries}, 
\ref{prop_grouSer}, and \ref{prop_prodKser}, numericly
$$
F_{B,\,h}(x)=
F_{B,\,v,\,h}(x)+F_{C,\,v,\,h}(x)^2F_{B,\,h}(x)\,\text{ for every }\,x\in(-1,\,1)\,.
$$
Claim~1 follows at once. If $F_{B,h}(x)\ne 0$, we can divide by it and get
the stated bound
$$
{\textstyle
\frac{F_{B,\,v,\,h}(x)}{F_{B,\,h}(x)}=1-F_{C,\,v,\,h}(x)^2\le1\,.
}
$$
\eproof

\subsection{General weights}\label{subsec_vne1}

In this subsection, we drop the convexity condition for weights.
Recall the conditional, absolute, and infinite sums 
$U(1)_c$, $U(1)_a$, and 
$U(1)_c=\pm\infty$ in Definitions~\ref{def_U12} and 
\ref{def_U1infty}. The generating functions $F_{A,v,h}(x)$, $F_{B,v,h}(x)$, 
$F_{C,v,h}(x)$, and $F_{E,v,h}(x)$ 
were defined in the previous subsection. Recall that all weights $h$ are $v$-transitive 
and that always $v\in V(h)$. In the next two theorems, 
the sums are conditional. In the former theorem, we assume the existence of the conditional 
$A_v$-, $B$-, $B_v$-, $C_v$-, and $D$-sums, and that $F_{B,h}(1)_c\ne0$. In the 
latter theorem, we replace this non-vanishing with the 
existence of the conditional $E_v$-sum. For $w,z\in\C$, we define the set
$$
w\cdot\mathrm{sqrt}(z)=\{w\al\cc\;\al\in\C,\,\al^2=z\}\ \ (\sus\C)\,.
$$
If $w=0$ or $z=0$, the set equals $\{0\}$. Else, it has two non-zero 
elements differing by sign. If $w=1$, we write just $\mathrm{sqrt}(z)$.

\begin{thm}[$U(1)_c$]\label{thm_vne1gen1}
Let $v\in\N\setminus\{1\}$ and let $h\cc\N_2\to\C$ be 
a~$v$-transitive light 
weight such that $v\in V(h)$. If the conditional sums $F_{A,v,h}(1)_c$, $F_{B,h}(1)_c$, $F_{B,v,h}(1)_c$, $F_{C,v,h}(1)_c$, and $F_{D,h}(1)_c$ exist and if $F_{B,h}(1)_c\ne0$, then
$$
F_{A,\,v,\,h}(1)_c\in
F_{D,\,h}(1)_c\cdot
\mathrm{sqrt}\Big(1-
\frac{F_{B,\,v,\,h}(1)_c}{F_{B,\,h}(1)_c}\Big)\ \ (\sus\C)\,.
$$  
\end{thm}
\proof
Identity~1 in Proposition~\ref{prop_12relvne1} states that $F_{A,v,h}(x)=F_{C,v,h}(x)F_{D,h}(x)$.
By Proposition~\ref{prop_Abel1}, 
$F_{A,v,h}(1)_c=F_{C,v,h}(1)_cF_{D,h}(1)_c$. By identity~2 in Proposition~\ref{prop_12relvne1}, $F_{C,v,h}(x)^2=1-
\frac{F_{B,v,h}(x)}{F_{B,h}(x)}$. Thus by Proposition~\ref{prop_Abel3}, 
$${\textstyle
\big(F_{C,\,v,\,h}(1)_c\big)^2=1-\frac{F_{B,\,v,\,h}(1)_c}{F_{B,\,h}(1)_c}\,.
}
$$
The stated formula follows. 
\eproof

\noindent
As we mentioned, in the following variant of the previous 
theorem we replace the non-vanishing of $F_{B,h}(1)_c$ 
with the existence of 
$F_{E,v,h}(1)_c$.

\begin{thm}[$U(1)_c$]\label{thm_vne1gen1Car}
Let $v\in\N\setminus\{1\}$ and let $h\cc\N_2\to\C$ be 
a~$v$-transitive light 
weight such that $v\in V(h)$. If the conditional sums $F_{A,v,h}(1)_c$, $F_{B,h}(1)_c$, 
$F_{B,v,h}(1)_c$, $F_{C,v,h}(1)_c$, $F_{D,h}(1)_c$, and $F_{E,v,h}(1)_c$
exist, then
$$
F_{A,\,v,\,h}(1)_c\in
F_{D,\,h}(1)_c\cdot
\mathrm{sqrt}\Big(1-
\frac{F_{B,\,v,\,h}(1)_c}{F_{B,\,h}(1)_c}\Big)\ \ (\sus\C)
$$  
provided that $F_{B,h}(1)_c\ne0$. If $F_{B,h}(1)_c=0$ then 
$F_{B,v,h}(1)_c=0$, and we interpret the expression $\frac{0}{0}$ as $1$, so that then $F_{A,v,h}(1)_c=0$.
\end{thm}
\proof
If $F_{B,h}(1)_c\ne0$, we use the previous proof. If $F_{B,h}(1)_c=0$, then it follows 
from identity~2 in Proposition~\ref{prop_12relvne1}, and 
from Propositions~\ref{prop_Abel1} and \ref{prop_totezPlus} that $F_{B,v,h}(1)_c=0$ as well. Identity~3 in Proposition~\ref{prop_12relvne1} 
and Proposition~\ref{prop_Abel1} yield 
$F_{C,v,h}(1)_c=0$. Identity~1 in Proposition~\ref{prop_12relvne1} and 
Proposition~\ref{prop_Abel1} give 
$F_{A,v,h}(1)_c=0$.
\eproof

Next, we obtain a~version of the previous results for absolute sums. It 
suffices to assume the existence of absolute $B$-, $B_v$-, and $D$-sums.

\begin{thm}[$U(1)_a$]\label{thm_vne1gen2}
Let $v\in\N\setminus\{1\}$ and let $h\cc\N_2\to\C$ be 
a~$v$-transitive light 
weight such that $v\in V(h)$. 
The following holds.
\begin{enumerate}
\item If $F_{B,v,h}(x)
=F_{B,h}(x)$ then $F_{A,v,h}(x)=0$.
\item Let $c:=\sum_{n\ge1}|h(B_n)|<1$, let $F_{B,v,h}(x)\ne F_{B,h}(x)$, let $k\in\N_0$ be minimum with 
$h(B_{v,k})\ne h(B_k)$, and let $t:=|h(B_{v,k})-h(B_k)|$. If the 
absolute sums $F_{B,v,h}(1)_a$ and 
$F_{D,h}(1)_a$ exist, and if
$$
{\textstyle \frac{1}{(1-c)t}\cdot
\big(ct+\sum_{n>k}|h(B_{v,\,n})-h(B_n)|\big)<1\,,
}
$$ 
then the absolute sum $F_{B,h}(1)_a$ exists, is nonzero, and
$$
F_{A,\,v,\,h}(1)_a\in
F_{D,\,h}(1)_a\cdot
\mathrm{sqrt}\Big(1-
\frac{F_{B,\,v,\,h}(1)_a}{F_{B,\,h}(1)_a}\Big)\ \ (\sus\C)\,.
$$  
\end{enumerate}
\end{thm}
\proof
1. If $F_{B,v,h}(x)=F_{B,h}(x)$, 
then by identity~2 in Proposition~\ref{prop_12relvne1} we have $F_{C,v,h}(x)^2=1-
\frac{F_{B,v,h}(x)}{F_{B,h}(x)}=1-1=0$ and $F_{C,v,h}(x)=0$. Identity~1 in Proposition~\ref{prop_12relvne1} gives $F_{A,v,h}(x)=
F_{C,v,h}(x)F_{D,h}(x)=0$. 

2. Let the stated hypothesis hold. Identity~2 in 
Proposition~\ref{prop_12relvne1} states that $F_{C,v,h}(x)^2=1-
\frac{F_{B,v,h}(x)}{F_{B,h}(x)}$. By
Proposition~\ref{prop_bezAbel3}, $F_{B,h}(1)_a\ne0$ and
$${\textstyle
\big(F_{C,\,v,\,h}(1)_a\big)^2=1-\frac{F_{B,\,v,\,h}(1)_a}{F_{B,\,h}(1)_a}\,. 
}
$$
Identity~1 in 
Proposition~\ref{prop_12relvne1} states that $F_{A,v,h}(x)=
F_{C,v,h}(x)F_{D,h}(x)$. By part~1 of  Proposition~\ref{prop_bezAbel1},  
$F_{A,v,h}(1)_a=F_{C,v,h}(1)_a
F_{D,h}(1)_a$. The stated formula follows. 
\eproof

In the last theorem in this subsection, we assume that the weights are complex nonnegative and that the 
conditional $D$-sum is $+\infty$. Now, as in the previous section, we 
also allow starting vertices $u\ne1$. 

\begin{thm}[$U(1)_c=+\infty$]\label{thm_vne1genInf}
Let $\xi\in\C_1$,  
$v\in\N\setminus\{1\}$, let  $h\cc\N_2\to\xi\R_{\ge0}$ be 
a~$v$-transitive light weight 
such that $v\in V(h)$, and let $(F_{D,h})_{1/\xi}(1)_c=+\infty$. The 
following holds.
\begin{enumerate}
\item Then also
$$
(F_{A,\,v,\,h})_{1/\xi}(1)_c=+\infty\,.
$$
\item If $u\in\N\setminus\{1,v\}$ is a~vertex such that $\{1,u\}$ is cycle-covered and $h$ is $u$-transitive, then also 
$$
(F_{A',\,v,\,h})_{1/\xi}(1)_c=+\infty\,.
$$
Recall that $F_{A',v,\,h}(x)=\sum_{n\ge0}h(A_{v,n}')x^n$ with $A_{v,n}'$ being the set of walks with length $n$ in $\vec{K}_{\N}$ 
that start at $u$ and visit $v$.
\item If the weight $h$ is in addition symmetric, then the 
assumption of cycle-cover in part~2 may 
be stated as the assumption that $u\in V(h)\setminus\{1,v\}$.
\item Without the 
assumption of cycle-cover, part~2 generally does not hold for general weights $h$.
\end{enumerate}
\end{thm}
\proof
1. Since $v\in V(h)$, there exists a~shortest $1$-$v$ walk $w_0$ with $h(w_0)\ne0$. It follows that $w_0\in 
C_{v,m}$, where $m=|w_0|$. Let $c=\xi^{-m}h(w_0)>0$. By the subseries
bound, $\xi^{-m}h(C_{v,m})\ge c>0$. Identity~1 in 
Proposition~\ref{prop_12relvne1} states that $h(A_{v,n})=\sum_{j=0}^n h(C_{v,j})h(D_{n-j})$,  which implies that
$$
\xi^{-n-m}h(A_{v,\,n+m})\ge 
c\cdot\xi^{-n}h(D_n)
$$ 
for every $n\in\N_0$. Hence 
$$
{\textstyle
\sum_{j=0}^n\xi^{-j-m}h(A_{v,\,j+m})
\ge c\cdot\sum_{j=0}^n\xi^{-j}
h(D_j)\to+\infty
}
$$
for $n\to\infty$. We see that $(F_{A,v,h})_{1/\xi}(1)_c=+\infty$. 

2. Since $\{1,u\}$ is cycle-covered, there exists 
a~$1$-$u$ walk $w_0$ and a~$u$-$1$ walk $w_1$ in $\vec{K}_{\N}$ such 
that $h(w_0)\ne0\ne h(w_1)$. For $n\in\N_0$, let $D_n'$ be the set 
of walks of length $n$ in $\vec{K}_{\N}$ starting at $u$. Due 
to the walk $w_0$ and Proposition~\ref{prop_moveStVer}, 
every series $\sum_{w\in D_n'}h(w)$, $n\in\N_0$, 
absolutely converges. The generating function 
$${\textstyle
F_{D',\,h}(x)=
\sum_{n\ge0}h(D_n')x^n
}
$$ 
is therefore well defined. We show that 
$(F_{D',h})_{1/\xi}(1)_c=+\infty$. Let $m=|w_1|$ ($\in\N_0$) and 
$c=\xi^{-m}h(w_1)>0$. For every walk $w\in D_n$, we consider the concatenated walk $w'=w_1w\in 
D_{m+n}'$. The subseries bound and Proposition~\ref{prop_LKseries} yield 
the lower bound
$$
\xi^{-m-n}h(D_{m+n}')\ge c\cdot 
\xi^{-n}h(D_n)\,.
$$
Thus, since $(F_{D,h})_{1/\xi}(1)_c=+\infty$, also $(F_{D',h})_{1/\xi}(1)_c=+\infty$.

By Proposition~\ref{prop_12relvne1}, the formal identity
$$
F_{A',\,v,\,h}(x)=F_{C',\,v,\,h}(x)F_{D',\,h}(x)
$$
holds. Since $v\in V(h)$ and we have the walk $w_1$, there exists 
a~$u$-$v$ walk with nonzero weight. 
Also, $(F_{D',h})_{1/\xi}(1)_c=+\infty$. We are therefore basically in the situation of part~1, only the starting vertex~$1$ is replaced with $u$, and we deduce that 
$$
(F_{A',\,v,\,h})_{1/\xi}(1)_c=+\infty\,.
$$

3. If $h$ is symmetric, then any nonzero $1$-$u$ walk yields a~nonzero 
cycle containing $\{1,u\}$. 

4. The counterexample has the weight $h$ defined in the proof of part~4 of 
Theorem~\ref{thm_v1conNN}. Then $F_{D',h}(1)_c=+\infty$, 
but $F_{A',v,h}(1)_c=F_{A',v,h}(1)_a=0$ for every $u,v\in\N\setminus\{1\}$ with $u>v$.
\eproof

\noindent
In parts~2--4 we assume that $u\ne v$, because the case $u=v$ was 
resolved in the previous section.

\subsection{Convex weights}\label{subsec_vne1C}

In this subsection, the convexity of weights is generalized to $\xi$-convexity. In the next theorem, we 
generalize the formula (\ref{eq_odmoc}) 
to complex nonnegative weights. 
We allow starting vertices $u\ne1$.
In fact, we present in a~general form the argument not given in 
Subsection~\ref{subsec_PTvSMvn0}.
Unlike in Theorem~\ref{thm_v1conNN}, we could not join the cases 
$(F_{B,h})_{1/\xi}(1)_c<+\infty$ and $(F_{B,h})_{1/\xi}(1)_c=+\infty$, 
because now the numerator is  $(F_{B,v,h})_{1/\xi}(1)_c$  and not 
$1$. We relegate the latter case to part~1 of 
Theorem~\ref{thm_vne1conInf}. 

\begin{thm}[$U(1)_c$ and  $\xi^{-1}h\ge0$]\label{thm_vne1conNN}
Suppose that $\xi\in\C_1$, that $v\in\N\setminus\{1\}$, that 
$h\cc\N_2\to\xi\R_{\ge0}$ is a~$v$-transitive $\xi$-convex light weight 
with $v\in V(h)$, and that $(F_{B,h})_{1/\xi}(1)_c<+\infty$. The following holds.
\begin{enumerate}
\item Then $(F_{B,h})_{1/\xi}(1)_c\ge1$, the conditional sum $(F_{B,v,h})_{1/\xi}(1)_c$ exists and 
$(F_{B,v,h})_{1/\xi}(1)_c\le (F_{B,h})_{1/\xi}(1)_c$, and 
$$
\lim_{n\to\infty}\xi^{-n}h(A_{v,\,n})=
\sqrt{1-\frac{(F_{B,\,v,\,h})_{1/\xi}(1)_c}{(F_{B,\,h})_{1/\xi}(1)_c}}\in[0,\,1)\,.
$$
\item  If $\{1,v\}$ is cycle-covered, which is the case when $h$
is symmetric, then the last interval may be replaced with $(0,1)$.
\item  If $u\in\N\setminus\{1,v\}$ is 
a~vertex such that $\{1,u\}$ is cycle-covered, then we have $(F_{B',h})_{1/\xi}(1)_c<+\infty$ and part~1 holds with 
$B$ replaced by $B'$ and $A_{v,n}$ by $A_{v,n}'$. Recall that the sets of 
walks $A_{v,n}'$ and the generating functions $F_{B',h}(x)$ and 
$F_{B',v,h}(x)$ are defined relative to the starting vertex $u$. If 
$\{1,v\}$ is cycle-covered, which is the case when $h$
is symmetric, then the displayed interval $[0,1)$ may be replaced with $(0,1)$.
\end{enumerate}
\end{thm}
\proof
1. $(F_{B,h})_{1/\xi}(1)_c\ge1$ because $h(B_0)=1$. 
The sum $(F_{B,v,h})_{1/\xi}(1)_c$ exists and 
$(F_{B,v,h})_{1/\xi}(1)_c\le (F_{B,h})_{1/\xi}(1)_c$ due to the 
subseries bound, because for every $n\in\N_0$, the series $S_n=\sum_{w\in B_{v,n}}\xi^{-n}h(w)$ is a~subseries of the series 
$S_n'=\sum_{w\in B_n}\xi^{-n}h(w)$ and both series have nonnegative real 
summands. We assume that $(F_{B,h})_{1/\xi}(1)_c<+\infty$ and 
prove the stated formula. Identity~2 in Proposition~\ref{prop_12relvne1} 
says that $F_{C,v,h}(x)^2=1-
\frac{F_{B,v,h}(x)}{F_{B,h}(x)}$. Equivalently and formally, 
$$
F_{C,\,v,\,h}(x)^2F_{B,\,h}(x)=F_{B,\,h}(x)-F_{B,\,v,\,h}(x)\,.
$$
By Proposition~\ref{prop_onD}, $F_{D,h}(x)=\frac{1}{1-\xi x}$. Since 
$B_n,B_{v,n},C_{v,n}\sus D_n$, the 
subseries bound implies that all coefficients of the generating 
functions $(F_{B,h})_{1/\xi}(x)$, $(F_{B,v,h})_{1/\xi}(x)$, and 
$(F_{C,v,h})_{1/\xi}(x)$ lie in $[0,1]$. These generating functions 
are therefore in $H_1$. By Propositions~\ref{prop_LKseries}, 
\ref{prop_grouSer}, and \ref{prop_proSer}, the numeric 
identity
$$
(F_{C,\,v,\,h})_{1/\xi}(x)^2(F_{B,\,h})_{1/\xi}(x)=(F_{B,\,h})_{1/\xi}(x)-(F_{B,\,v,\,h})_{1/\xi}(x)
$$
holds for every number $x\in[0,1)$.
As we know, $(F_{B,h})_{1/\xi}(x)\ge1$ for 
every such~$x$. We divide by  
$(F_{B,h})_{1/\xi}(x)$ and get the main numeric identity that
$$
{\textstyle
(F_{C,\,v,\,h})_{1/\xi}(x)^2=
1-\frac{(F_{B,\,v,\,h})_{1/\xi}(x)}{(F_{B,\,h})_{1/\xi}(x)}
\,\text{ for every }\,x\in[0,\,1)\,.
}
$$
We compute that
\begin{eqnarray*}
\big((F_{C,\,v,\,h})_{1/\xi}(1)_c\big)^2&=&
((F_{C,\,v,\,h})_{1/\xi})^2(1)_c=
\lim_{x\nearrow1}(F_{C,\,v,\,h})_{1/\xi}(x)^2\\
&=&\lim_{x\nearrow1}
{\textstyle
\big(1-\frac{(F_{B,\,v,\,h})_{1/\xi}(x)}{(F_{B,\,h})_{1/\xi}(x)}\big)=
1-
\frac{\lim_{x\nearrow1}(F_{B,\,v,\,h})_{1/\xi}(x)}{\lim_{x\nearrow1}(F_{B,\,h})_{1/\xi}(x)}}\\
&=&{\textstyle
1-\frac{(F_{B,\,v,h})_{1/\xi}(1)_c}{(F_{B,\,h})_{1/\xi}(1)_c}\,.
}
\end{eqnarray*}
The first equality follows from Proposition~\ref{pridat}. The second 
equality follows from Theorem~\ref{thm_abel3}. The third 
equality follows from the above main numeric identity. The fifth and 
fourth equalities follow from Theorem~\ref{thm_abel3} and the 
arithmetic of functional limits. Since 
$(F_{C,v,\,h})_{1/\xi}(1)_c\ge0$, the number $(F_{C,v,\,h})_{1/\xi}(1)_c$ is equal to the displayed root. By identity~1 in Proposition~\ref{prop_12relvne1} 
and by  
Proposition~\ref{prop_onD}, $F_{A,v,h}(x)=F_{C,v,h}(x)F_{D,h}(x)=F_{C,v,h}(x)\frac{1}{1-\xi x}$. By Proposition~\ref{prop_FPS1},
$$
\lim_{n\to\infty}\xi^{-n}h(A_{v,\,n})=
(F_{C,\,v,\,h})_{1/\xi}(1)_c
$$ 
and the stated formula follows. 

2. We consider the series $S_n$ and $S_n'$ defined at the beginning of 
the proof of part~1. As we know, they have nonnegative summands and the 
former is a~subseries of the latter. Since $\{1,v\}$ is cycle-covered, for some $n$, the latter series contains a~positive summand 
that is not present in the former 
series. Thus, $(F_{B,v,h})_{1/\xi}(1)_c<(F_{B,h})_{1/\xi}(1)_c$ and
the root is positive. For symmetric~$h$, any nonzero $1$-$v$ 
walk yields a~nonzero cycle containing $1$ and $v$.  

3. Let $u$ be as stated. Since there is a~nonzero $1$-$u$ walk $w_0$, by 
Proposition~\ref{prop_moveStVer} the coefficients of $F_{B',h}(x)$ are well 
defined. We show that $(F_{B',h})_{1/\xi}(1)_c<+\infty$. This follows from the assumed bound   
$(F_{B,h})_{1/\xi}(1)_c<+\infty$, the subseries bound, and from the existence of $w_0$ and a~nonzero $u$-$1$ walk $w_1$, because 
every $u$-$u$ walk $w$ can be extended to a~$1$-$1$ walk $w_0ww_1$. 
It is easy to see that the three initial claims in part~1 also hold 
with $B$ replaced by $B'$, and we get the formula of part~1 relative to the 
starting vertex $u$. Finally, if $\{1,v\}$ is cycle-covered, then so is 
$\{u,v\}$ because $\{1,u\}$ is cycle-covered.   
\eproof

The next two theorems treat conditional sums. The 
weights are complex and $\xi$-convex, and we generalize the formula 
(\ref{eq_odmoc}). The formula (\ref{eq_exactly1}) is generalized in 
Theorem~\ref{thm_vne1conInf}. As before, $0$-convex weights are 
treated separately. In the former theorem, we assume the existence of  
conditional 
$B_v$-, $B$-, and $C_v$-sums, and that $(F_{B,h})_{1/\xi}(1)_c\ne0$. In the 
latter theorem, we replace this non-vanishing with the existence of 
conditional $E_v$-sum. 

\begin{thm}[$U(1)_c$]\label{thm_vne1con1}
Let $v\in\N\setminus\{1\}$ and let 
$h\cc\N_2\to\C$ be a~$v$-transitive 
$\xi$-convex light weight such that $v\in V(h)$. The following holds.
\begin{enumerate}
\item 
If $\xi\in\C_1$ and the conditional sums 
$(F_{B,v,h})_{1/\xi}(1)_c$,  $(F_{B,h})_{1/\xi}(1)_c$, and $(F_{C,v,h})_{1/\xi}(1)_c$ exist, and if $(F_{B,h})_{1/\xi}(1)_c\ne0$, then
$$
\lim_{n\to\infty}
\xi^{-n}h(A_{v,\,n})
\in\mathrm{sqrt}\Big(1-\frac{(F_{B,\,v,\,h})_{1/\xi}(1)_c}{(F_{B,\,h})_{1/\xi}(1)_c}\Big)\ \ (\sus\C)\,.
$$
\item If $\xi=0$ and the conditional sums $F_{B,v,h}(1)_c$, 
$F_{B,h}(1)_c$, and $F_{C,v,h}(1)_c$ exist, and if $F_{B,h}(1)_c\ne0$, then the conditional sum $F_{A,v,h}(1)_c$ exists and
$$
F_{A,\,v,\,h}(1)_c
\in\mathrm{sqrt}\Big(1-\frac{F_{B,\,v,\,h}(1)_c}{F_{B,\,h}(1)_c}\Big)\ \ (\sus\C)\,.
$$
\end{enumerate}
\end{thm}
\proof
1. By identity~1 in Proposition~\ref{prop_12relvne1} and by Proposition~\ref{prop_onD}, 
we have $F_{A,v,h}(x)=F_{C,v,h}(x)F_{D,h}(x)=F_{C,v,h}(x)\frac{1}{1-\xi x}$. By Proposition~\ref{prop_FPS1},
$$
\lim_{n\to\infty}\xi^{-n}h(A_{v,\,n})=
(F_{C,\,v,\,h})_{1/\xi}(1)_c\,.
$$ 
Identity~2 in Proposition~\ref{prop_12relvne1} says that $(F_{C,v,h})_{1/\xi}(x)^2=1-
\frac{(F_{B,v,h})_{1/\xi}(x)}{(F_{B,h})_{1/\xi}(x)}$. By
Proposition~\ref{prop_Abel3},  
$${\textstyle
\big((F_{C,\,v,\,h})_{1/\xi}(1)_c\big)^2=1-\frac{(F_{B,\,v,\,h})_{1/\xi}(1)_c}{(F_{B,\,h})_{1/\xi}(1)_c}\,. 
}
$$
The stated formula follows.

2. By identity~1 in 
Proposition~\ref{prop_12relvne1} and by Proposition~\ref{prop_onD},   
$F_{A,v,h}(x)=F_{C,v,h}(x)$. We use the displayed formula for
$(F_{C,v,h})_{1/\xi}(1)_c$ with $\xi=1$.
\eproof

\noindent
As we mentioned, we replace the assumption that $F_{B,h}(1)_c\ne0$ 
with the assumption that the conditional sum   
$F_{E,v,h}(1)_c$ exists and obtain the following variant of 
Theorem~\ref{thm_vne1con1}.

\begin{thm}[$U(1)_c$]\label{thm_vne1con1Car}
Let $v\in\N\setminus\{1\}$ and let  
$h\cc\N_2\to\C$ be a~$v$-transitive 
$\xi$-convex light weight such that $v\in V(h)$. The following holds.
\begin{enumerate}
\item If $\xi\in\C_1$ and if the four conditional sums $(F_{B,v,h})_{1/\xi}(1)_c$, $(F_{B,h})_{1/\xi}(1)_c$, $(F_{C,v,h})_{1/\xi}(1)_c$, and $(F_{E,v,h})_{1/\xi}(1)_c$ exist, then
$$
\lim_{n\to\infty}
\xi^{-n}h(A_{v,\,n})
\in\mathrm{sqrt}\Big(1-\frac{(F_{B,\,v,\,h})_{1/\xi}(1)_c}{(F_{B,\,h})_{1/\xi}(1)_c}\Big)\ \ (\sus\C)
$$ 
if $(F_{B,h})_{1/\xi}(1)_c\ne0$. For $(F_{B,h})_{1/\xi}(1)_c=0$ we have 
$(F_{B,v,h})_{1/\xi}(1)_c=0$ and interpret the expression 
$\frac{0}{0}$ as $1$, so that the limit is zero.
\item If $\xi=0$ and the conditional sums $B_{v,h}(1)_c$, $B_h(1)_c$, $C_{v,h}(1)_c$, and $E_{v,h}(1)_c$ exist, then the conditional sum $F_{A,v,h}(1)_c$ exists and
$$
F_{A,\,v,\,h}(1)_c
\in\mathrm{sqrt}\Big(1-\frac{F_{B,\,v,\,h}(1)_c}{F_{B,\,h}(1)_c}\Big)\ \ (\sus\C)
$$ 
if $F_{B,h}(1)_c\ne0$. For $F_{B,h}(1)_c=0$ we have 
$F_{B,v,h}(1)_c=0$ and interpret the expression $\frac{0}{0}$ as $1$, so that $F_{A,v,h}(1)_c=0$.
\end{enumerate}
\end{thm}
\proof
If $(F_{B,h})_{1/\xi}(1)_c\ne0$, respectively $F_{B,h}(1)_c\ne0$, we 
are in the previous theorem. We therefore assume that $(F_{B,h})_{1/\xi}(1)_c=0$, 
respectively that $F_{B,h}(1)_c=0$.

1. From the proof of part~1 of Theorem~\ref{thm_vne1con1} we know that 
$$
\lim_{n\to\infty}\xi^{-n}h(A_{v,\,n})=
(F_{C,\,v,\,h})_{1/\xi}(1)_c\,.
$$ 
Using identity~2 in Proposition~\ref{prop_12relvne1} and 
Propositions~\ref{prop_Abel1} and \ref{prop_totezPlus}, we deduce from $(F_{B,h})_{1/\xi}(1)_c=0$ that $(F_{B,v,h})_{1/\xi}(1)_c=0$. Identity~3 in Proposition~\ref{prop_12relvne1} 
and Proposition~\ref{prop_Abel1} give that 
$(F_{C,v,h})_{1/\xi}(1)_c=0$. The displayed limit is therefore $0$.

2. Since $\xi=0$, identity~1 in 
Proposition~\ref{prop_12relvne1} and Proposition~\ref{prop_onD} yield 
$F_{A,v,h}(x)=F_{C,v,h}(x)$. As we know, $F_{B,h}(1)_c=0$ 
implies that $F_{B,v,h}(1)_c=
F_{C,v,h}(1)_c=0$. So $F_{A,v,h}(1)_c=0$.
\eproof

Next, we obtain a~version of the previous results for absolute sums. 
It suffices to assume the existence of absolute $B$- and $B_v$-sums.

\begin{thm}[$U(1)_a$]\label{thm_vne1con2}
Let $v\in\N\setminus\{1\}$ and let $h\cc\N_2\to\C$ be a~$v$-transitive $\xi$-convex light 
weight such that $v\in V(h)$. The following holds.
\begin{enumerate}
\item If $F_{B,v,h}(x)=F_{B,h}(x)$, then $F_{A,v,h}(x)=0$.
\item Suppose that $c:=\sum_{n\ge1}|h(B_n)|<1$, that 
$F_{B,v,h}(x)\ne F_{B,h}(x)$, that $k$ in $\N_0$ is minimum with 
$h(B_{v,k})\ne h(B_k)$, and that $t:=|h(B_{v,k})-h(B_k)|$.
\begin{enumerate}
\item If $\xi\in\C_1$, if the absolute sum $(F_{B,v,h})_{1/\xi}(1)_a$ exists, and if 
\begin{equation*}
{\textstyle
\frac{1}{(1-c)t}
\cdot\big(ct+\sum_{n>k}|h(B_{v,\,n})-h(B_n)|\big)<1\,,\tag{$*$}
}
\end{equation*} 
then the absolute sum $(F_{B,h})_{1/\xi}(1)_a$ exists, is nonzero, and
$$
\lim_{n\to\infty}\xi^{-n}h(A_{v,\,n})\in
\mathrm{sqrt}\Big(1-
\frac{(F_{B,\,v,\,h})_{1/\xi}(1)_a}{(F_{B,\,h})_{1/\xi}(1)_a}\Big)\ \ (\sus\C)\,.
$$
\item If $\xi=0$, if the absolute sum $F_{B,v,h}(1)_a$ exists, and if  the bound $(*)$ holds,
then the absolute sum $F_{B,h}(1)_a$ exists, is nonzero, and
$$
F_{A,\,v,\,h}(1)_a\in
\mathrm{sqrt}\Big(1-
\frac{F_{B,\,v,\,h}(1)_a}{F_{B,\,h}(1)_a}\Big)\ \ (\sus\C)\,.
$$
\end{enumerate}
\end{enumerate}
\end{thm}
\proof
1. Let $F_{B,v,h}(x)=F_{B,h}(x)$. 
By identity~2 in Proposition~\ref{prop_12relvne1} we 
have $F_{C,v,h}(x)^2=1-
\frac{F_{B,v,h}(x)}{F_{B,h}(x)}=
1-1=0$ and $F_{C,v,h}(x)=0$. By identity~1 in 
Proposition~\ref{prop_12relvne1},   $F_{A,v,h}(x)=
F_{C,v,h}(x)F_{D,h}(x)=0$.

2 (a). Let the stated hypothesis hold. By identity~2 in 
Proposition~\ref{prop_12relvne1}, formally  $F_{C,v,h}(x)^2=1-
\frac{F_{B,v,h}(x)}{F_{B,h}(x)}$. By Proposition~\ref{prop_bezAbel3}, the sum $(F_{B,h})_{1/\xi}(1)_a$ exists,  is nonzero, and
$${\textstyle
\big((F_{C,v,h})_{1/\xi}(1)_a\big)^2=1-\frac{(F_{B,v,h})_{1/\xi}(1)_a}{(F_{B,h})_{1/\xi}(1)_a}\,.
}
$$
By identity~1 in 
Proposition~\ref{prop_12relvne1} and by  
Proposition~\ref{prop_onD}, we have 
$F_{A,v,h}(x)=F_{C,v,h}(x)F_{D,h}(x)=F_{C,v,h}(x)\frac{1}{1-\xi x}$.
Proposition~\ref{prop_FPS1} implies the stated formula: $\lim_{n\to\infty}\xi^{-n}h(A_{v,n})$ equals
$$
{\textstyle
(F_{C,\,v,\,h})_{1/\xi}(1)_c=(F_{C,\,v,\,h})_{1/\xi}(1)_a
\in\mathrm{sqrt}\big(1-\frac{(F_{B,\,v,\,h})_{1/\xi}(1)_a}{(F_{B,\,h})_{1/\xi}(1)_a}\big)\,.
}
$$

2 (b). Let the stated hypothesis hold. By part~1 used for $\xi=1$, we 
have $\big(F_{C,v,h}(1)_a\big)^2=1-\frac{F_{B,v,h}(1)_a}{F_{B,h}(1)_a}$. 
Identity~1 in 
Proposition~\ref{prop_12relvne1} and Proposition~\ref{prop_onD} imply that 
$F_{A,v,h}(x)=F_{C,v,h}(x)F_{D,h}(x)=F_{C,v,h}(x)$. We see that the stated formula holds.
\eproof

The last theorem of this article treats the case $U(1)_c=\pm\infty$.
We allow starting vertices $u\ne1$.

\begin{thm}[$U(1)_c=\pm\infty$]\label{thm_vne1conInf}
Let $v\in\N\setminus\{1\}$, let
$\xi\in\C_1$, and let $h\cc\N_2\to\xi\R$ be a~$v$-transitive light weight such that $v\in V(h)$. The following holds.
\begin{enumerate}
\item If $h\cc\N_2\to\xi\R_{\ge0}$, if $h$ is $\xi$-convex, and if $(F_{B,h})_{1/\xi}(1)_c=+\infty$, then
$$
\lim_{n\to\infty}\xi^{-n}h(A_{v,\,n})=1\,.
$$
\item Under the assumptions of part $1$ on $h$, if $u\in\N\setminus\{1,v\}$ is a~vertex such that 
$\{1,u\}$ is cycle-covered and that $h$ is $u$-transitive, then also
$$
\lim_{n\to\infty}\xi^{-n}h(A_{v,\,n}')=1\,.
$$
Recall that $A_{v,n}'$ is the set of walks with length $n$ in $\vec{K}_{\N}$ that 
start at~$u$ and visit $v$.
If $h$ is in addition symmetric, it suffices to assume that $u$ is in 
$V(h)\setminus\{1,v\}$.
\item If $h$ is $0$-convex, if $F_{B,h}(x)$ and $F_{C,v,h}(x)$ are in $H_1$, if the conditional sum $(F_{B,v,h})_{1/\xi}(1)_c$ exists and is finite, and if $(F_{B,h})_{1/\xi}(1)_c=\pm\infty$, then  
$$
\lim_{x\nearrow1}(F_{A,\,v,\,h})_{1/\xi}(x)\in\{-1,\,1\}\,.
$$ 
\end{enumerate}
\end{thm}
\proof
1. By Proposition~\ref{prop_onD},  
$F_{D,h}(x)=\frac{1}{1-\xi x}$. 
The subseries bound implies that for every $n\in\N_0$,
$$
0\le\xi^{-n}h(B_n),\,
\xi^{-n}h(B_{v,\,n}),\,
\xi^{-n}h(C_{v,\,n}),\,\xi^{-n}h(E_{v,\,n})
\le\xi^{-n}h(D_n)=1\,.
$$ 
The generating functions 
$F_{B,h}(x)$, $F_{B,v,h}(x)$, $F_{C,v,h}(x)$, and $F_{E,v,h}(x)$ 
are therefore in $H_1$. Two facts are clear.
(i) By identity~2 in Proposition~\ref{prop_12relvne1}, the 
function $(F_{C,v,h})_{1/\xi}(x)$  is nonnegative and 
bounded for $x\in[0,1)$. (ii) The function 
$(F_{E,v,h})_{1/\xi}(x)$ is positive and increasing for 
$x\in(0,1)$: any shortest nonzero $1$-$v$ walk shows that 
$\xi^{-m}h(E_{v,m})>0$ for some 
$m\in\N$, which gives a~positive (non-constant) coefficient in $(F_{E,v,h})_{1/\xi}(x)$. Using (i), (ii), and identity~3 in 
Proposition~\ref{prop_12relvne1}, we deduce that the function 
$(F_{B,v,h})_{1/\xi}(x)$ is nonnegative and bounded for  $x\in[0,1)$. We compute that
\begin{eqnarray*}
\big((F_{C,\,v,\,h})_{1/\xi}(1)_c\big)^2
&=&((F_{C,\,v,\,h})_{1/\xi})^2(1)_c=
\lim_{x\nearrow1}(F_{C,\,v,\,h})_{1/\xi}(x)^2\\
&=&\lim_{x\nearrow1}
{\textstyle
\big(1-\frac{(F_{B,\,v,\,h})_{1/\xi}(x)}{(F_{B,\,h})_{1/\xi}(x)}\big)=1-
\frac{\lim_{x\nearrow1}(F_{B,\,v,\,h})_{1/\xi}(x)}{\lim_{x\nearrow1}(F_{B,\,h})_{1/\xi}(x)}}\\
&=&{\textstyle
1-\frac{c}{(F_{B,\,h})_{1/\xi}(1)_c}=1-\frac{c}{+\infty}=1\,.
}
\end{eqnarray*}
The first equality follows from Proposition~\ref{pridat}. The second 
equality follows from Theorem~\ref{thm_abel3}. The third 
equality follows from the main numeric identity obtained in the 
proof of Theorem~\ref{thm_vne1conNN}. The fifth and fourth equalities follow from the boundedness of the function 
$(F_{B,v,h})_{1/\xi}(x)$, from Theorem~\ref{thm_abel3}, and from the 
arithmetic of functional limits. The 
sixth equality follows from the assumptions. The last, seventh 
equality follows from the arithmetic of infinities. Since 
$(F_{C,v,\,h})_{1/\xi}(1)_c\ge0$, the number $(F_{C,v,\,h})_{1/\xi}(1)_c=1$. By identity~1 in 
Proposition~\ref{prop_12relvne1} and by Proposition~\ref{prop_onD}, we 
have equalities  
$F_{A,v,h}(x)=F_{C,v,h}(x)F_{D,h}(x)=
F_{C,v,h}(x)\frac{1}{1-\xi x}$. By 
Proposition~\ref{prop_FPS1}, 
$$
\lim_{n\to\infty}
\xi^{-n}h(A_{v,\,n})=(F_{C,\,v,\,h})_{1/\xi}(1)_c=1\,.
$$

2. Since $\{1,u\}$ is cycle-covered, there exists 
a~$1$-$u$ walk $w_0$ and a~$u$-$1$ walk $w_1$ such that $h(w_0)\ne0\ne 
h(w_1)$. Thus $u\in V(h)$ and the generating function $F_{B',h}(x)$ is 
well defined. We show that 
$(F_{B',h})_{1/\xi}(1)_c=+\infty$. This follows from the assumption that 
$(F_{B,h})_{1/\xi}(1)_c=+\infty$ and from the inequality that
$$
\xi^{-m_1-m_0-n}h(B_{m_1+m_0+n}')\ge
c\cdot\xi^{-n}h(B_n)
$$
for every $n\in\N_0$, where $c=\xi^{-m_1-m_0}h(w_1)h(w_0)>0$ and 
$m_i=|w_i|$. We obtain this inequality by considering, for every 
cycle $w$ covering $\{1\}$, the cycle $w'=w_1ww_0$ that covers $\{u\}$. 
Now we can apply part~1, if we show that the weight $h$ satisfies the 
modified hypothesis in which the starting vertex $1$ is replaced with $u$. As we 
know from the end of the proof of Proposition~\ref{prop_12relvne1}, 
there exists a~weight-preserving bijection $f_1\cc\N\to\N$ with 
$f_1(u)=v$. By the remark after Definition~\ref{def_conv}, 
$\xi$-convexity of $h$ is not affected by the change of the 
starting vertex because $u\in V(h)$. Finally, we 
have shown that 
$(F_{B',h})_{1/\xi}(1)_c=+\infty$. Thus, by part~1, the limit in part~2 
is $1$. If $h$ is symmetric and $u\in V(h)\setminus\{1\}$, then 
$\{1,u\}$ is cycle-covered by $w_0\overline{w_0}$. 

3. By identity~1 in Proposition~\ref{prop_12relvne1} and by Proposition~\ref{prop_onD}, 
$F_{A,v,h}(x)=F_{C,v,h}(x)$. By identity~2 in 
Proposition~\ref{prop_12relvne1},
$F_{C,v,h}(x)^2=1-\frac{F_{B,v,h}(x)}{F_{B,h}(x)}$. Equivalently and formally,
$$
(F_{C,\,v,\,h})_{1/\xi}(x)^2(F_{B,\,h})_{1/\xi}(x)=
(F_{B,\,h})_{1/\xi}(x)-(F_{B,\,v,\,h}(x))_{1/\xi}\,.
$$
By the assumption and by Proposition~\ref{prop_ConpowSer}, the generating functions $F_{B,h}(x)$, $F_{B,v,h}(x)$, and $F_{C,v,h}(x)$ are in $H_1$. Thus, by 
Propositions~\ref{prop_LKseries}, 
\ref{prop_grouSer}, and \ref{prop_proSer}, the displayed 
identity also holds numericly for every number $x\in[0,1)$. It follows 
from the assumption $(F_{B,h})_{1/\xi}(1)_c=\pm\infty$ that for every $x\in[0,1)$ close 
to $1$ we can divide by $(F_{B,h})_{1/\xi}(x)$. The main numeric identity
$$
{\textstyle
(F_{C,\,v,\,h})_{1/\xi}(x)^2=
1-\frac{(F_{B,\,v,\,h})_{1/\xi}(x)}{(F_{B,\,h})_{1/\xi}(x)}
}
$$
therefore holds for every $x\in(1-\de,1)$ and some $\de>0$.

We compute that
\begin{eqnarray*}
\lim_{x\nearrow1}(F_{A,\,v,\,h})_{1/\xi}(x)^2&=&    
\lim_{x\nearrow1}(F_{C,\,v,\,h})_{1/\xi}(x)^2=
\lim_{x\nearrow1}{\textstyle
\big(1-\frac{(F_{B,\,v,\,h})_{1/\xi}(x)}{(F_{B,\,h})_{1/\xi}(x)}\big)}\\
&=&{\textstyle
1-\frac{\lim_{x\nearrow1}(F_{B,\,v,\,h})_{1/\xi}(x)}{\lim_{x\nearrow1}(F_{B,\,h})_{1/\xi}(x)}=
1-\frac{(F_{B,\,v,\,h})_{1/\xi}(1)_c}{(F_{B,\,h})_{1/\xi}(1)_c}}\\
&=&{\textstyle
1-\frac{c}{\pm\infty}=1\,.
}
\end{eqnarray*}
We justified the first equality at the beginning of this proof.  
The second equality follows from the above main numeric identity. In the fourth and third equalities, we use
Theorem~\ref{thm_abel} in the numerator, Theorem~\ref{thm_abel2} 
in the denominator, and we use the arithmetic of functional limits. The fifth  equality follows from the 
assumptions. The last, sixth equality follows from the arithmetic of 
infinities. Since the function $(F_{A,v,h})_{1/\xi}(x)=(F_{C,v,h})_{1/\xi}(x)$ is continuous  for $x\in[0,1)$, we deduce that the finite limit 
$\lim_{x\nearrow1}(F_{A,v,h})_{1/\xi}(x)$ exists and
$$
\lim_{x\nearrow1}(F_{A,\,v,\,h})_{1/\xi}(x)=\lim_{x\nearrow1}(F_{C,\,v,\,h})_{1/\xi}(x)
\in\{-1,\,1\}\,.
$$   
\eproof  

\noindent
We are not completely satisfied with part~3 of the previous theorem and 
pose the following problem. 

\begin{prob}
May one omit in part~3 the assumption  of existence of the finite conditional sum  
$(B_{v,h})_{1/\xi}(1)_c$? That 
is, is it true that if $v\in\N\setminus\{1\}$, 
$\xi\in\C_1$, $h\cc\N_2\to\xi\R$ is a~$v$-transitive $0$-convex light weight such that $v\in V(h)$, and if  $(F_{B,h})_{1/\xi}(1)_c=\pm\infty$, then always
$$
\lim_{x\nearrow1}{\textstyle\frac{(F_{B,\,v,\,h})_{1/\xi}(x)}{(F_{B,\,h})_{1/\xi}(x)}=0\,?
}
$$    
\end{prob}

\section{Acknowledgments, AI, conclusion}\label{sec_concl}

The first author thanks C.~Banderier for helpful and inspiring discussions 
on this article. It was prepared with no  use of AI. To conclude and sum 
up, we demonstrated in our article that random walks can, in some 
situations, be defined and successfully investigated without 
any use of probability spaces. As we mentioned, in \cite{klaz2} we look into this 
matter more systematically.        

\appendix

\section{Sums $F_B(1)$ in P\'olya's theorem}\label{sec_sumsB}

In order to make the argument in Subsection~\ref{subsec_PTvSMv0} 
self-contained, we show 
in this appendix
that the sum $F_B(1)$ associated with the grid graph $G_d$ is $+\infty$ for $d\le2$, and lies in the interval $(1,+\infty)$ for $d\ge3$. Recall that 
$${\textstyle
F_B(1)=F_B^d(1)=\sum_{n\ge0}
\sum_{w\in B_n}(2d)^{-|w|}
=\sum_{n\ge0}(2d)^{-n}|B_n|\,,
}
$$
where $B_n=B_n^d$ is the set of walks with length $n$ in $G_d=\langle\Z^d,E_d\rangle$ that start 
and end at the origin $\overline{0}=\langle0,\ds,0\rangle$.

For the next proposition, recall that
for $n\in\N_0$ and $m_1,\ds.m_r\in\N_0$ such that $r\in\N$ and $m_1+\ds+m_r=n$, the multinomial coefficient
$$
{\textstyle
\binom{n}{m_1,\,m_2,\,\ds,\,m_r}=\frac{n!}{m_1!\cdot m_2!\cdot\ldots\cdot m_r!}
}
$$
counts the $r$-tuples $\langle X_1,\ds,X_r\rangle$ of mutually disjoint sets $X_i\sus\{1,\ds,n\}$
such that $|X_i|=m_i$ and $\bigcup_{i=1}^r X_i=\{1,\ds,n\}$.

\begin{prop}\label{prop_B_n}
Let $n\in\N_0$ and $d\in\N$. Then $|B_n^d|=0$ for every odd $n$, and for every even $n=2m$ we have
$$
|B_n^d|=|B_{2m}^d|=\sum_{\substack{j_1,\,\ds,\,j_d\in\N_0\\j_1+\ds+j_d=m}}
\binom{2m}{j_1,\,j_1,\,
j_2,\,j_2,\,\ds,\,j_d,\,j_d}\,.
$$
\end{prop}
\proof
Let $n\in\N_0$, $d\in\N$, and let $w\in B_n^d$. We write 
$w=\langle\overline{u}_0,\ds,
\overline{u}_n\rangle$, with $\overline{u}_i\in\Z^d$. For $k\in\{1,\ds,d\}$, we denote by 
$e_k=\langle0,\ds,0,1,0,\ds,0\rangle$
the $d$-tuple with $1$ on the $k$-th place and zeros elsewhere, and 
define
$$
j_k(w)=|\{i\in\{1,\,\ds,\,n\}\cc\;
\overline{u}_i-\overline{u}_{i-1}=e_k\}|\ \ (\in\N_0)\,.
$$
We similarly define $j_k'(w)$, by replacing $e_k$ with $-e_k$. It 
follows that $j_k(w)=j_k'(w)$ and $\sum_{k=1}^d(j_k(w)+j_k'(w))=n$ for 
every $k$ and every $w\in B_n^d$. Thus $n$ is even and $B_n^d=\emptyset$ for 
odd $n$. For $n=2m$ and nonnegative integers 
$j_1$, $\ds$, $j_d$ such that $\sum_{i=1}^d j_i=m$, the stated multinomial coefficient counts the 
walks $w\in B_n^d$ such that $j_k(w)=j_k'(w)=j_k$ for $k=1,\ds,d$. 
\eproof

For $d\le2$, we need a~lower bound on $|B_n^d|$. There are closed formulas for these cardinalities.

\begin{prop}\label{prop_dMost2}
Let $n=2m\in\N_0$ be even. Then
$$
|B_{2m}^1|=\binom{2m}{m}
\,\text{ and }\,|B_{2m}^2|=
\binom{2m}{m}^2\,.
$$
\end{prop}
\proof
Both formulas follow  from Proposition~\ref{prop_B_n}, the 
former immediately. As for the latter, Proposition~\ref{prop_B_n}
implies that
\begin{eqnarray*}
|B_{2m}^2|&=&{\textstyle
\sum_{j=0}^m\binom{2m}{j,\,j,\,m-j,\,m-j}=\binom{2m}{m}\sum_{j=0}^m
\binom{m}{j}^2\,.
}    
\end{eqnarray*}
By a~well known combinatorial argument, $\sum_{j=0}^m
\binom{m}{j}^2=\sum_{j=0}^m
\binom{m}{j}\binom{m}{m-j}=\binom{2m}{m}$.
\eproof

For $d\ge3$, such closed formulas are not known and probably do not exist, 
but good upper bounds fully suffice. We need three lemmas. The first 
one is an explicit form of a~well known asymptotics.

\begin{lem}[Stirling's formula]\label{lem_Stirling}
For every $n\in\N$, there exists $\theta_n\in[0,1]$ such that
$$
n!=1\cdot2\cdot\ldots\cdot n=
\exp(1/(12n+\theta_n))\sqrt{2\pi n}\big(n/\mathrm{e}\big)^n\,.
$$
\end{lem}
\proof
See \cite{robb}.
\eproof

\begin{lem}\label{lem_maxim}
Let $n\in\N_0$ and $r\in\N$. Then there exist numbers $n_1,\ds,n_r\in\N_0$ such that $n_1+\ds+n_r=n$ and $n_i\in\{\lfloor n/r\rfloor,\lceil n/r\rceil\}$
for every $i=1,\ds,r$. Moreover, 
$$
\binom{n}{m_1,\,\ds,\,m_r}\le
\binom{n}{n_1,\,\ds,\,n_r}
$$
for every $m_1,\ds,m_r\in\N_0$ 
with $m_1+\ds+m_r=n$.
\end{lem}
\proof
We prove both claims simultaneously and take any $m_i\in\N_0$, 
$i=1,\ds,r$, with $\sum_{i=1}^r m_i=n$. Suppose that there exist indices $j,j'\in
\{1,\ds,r\}$, $j\ne j'$, such that $m_j\le m_{j'}-2$. Let $m_i':=m_i$
for $i\in\{1,\ds,r\}\setminus\{j,j'\}$, $m_j':=m_j+1$, and 
$m_{j'}':=m_{j'}-1$. Then $\sum_{i=1}^r m_i'=n$ and
$$
{\textstyle
\binom{n}{m_1,\,\ds,\,m_r}<
\binom{n}{m_1',\,\ds,\,m_r'}
}
$$
because $\frac{1}{m_j!\cdot 
m_{j'}!}<\frac{1}{m_j'!\cdot m_{j'}'!}$ $\iff$ 
$\frac{m_j+1}{m_{j'}}<1$. 
Thus, the multinomial coefficient is maximized by an $r$-tuple of 
numbers $n_i\in\N_0$ satisfying $\sum_{i=1}^r n_i=n$ and $|n_j-
n_{j'}|\le1$ for every two indices $j,j'\in\{1,\ds,r\}$. So
there exists $s\in\N_0$ such that for every $i=1,\ds,r$, we have 
$n_i=s$ or $n_i=s+1$. We deduce that 
$n_i\in\{\lfloor n/r\rfloor,\lceil n/r\rceil\}$ for every $i=1,\ds,r$. 
\eproof

\begin{lem}\label{lem_treti}
Let $a_1,\ds,a_r\ge0$ be $r\in\N$ real numbers with $\sum_{i=1}^r a_i\le1$.  Then
$$
a_1^2+\ds+a_r^2\le\max_{1\le i\le r} a_i\,.
$$
\end{lem}
\proof
Indeed, 
$${\textstyle
\sum_{i=1}^r a_i^2\le
\sum_{i=1}^r(\max_{1\le i\le r} a_i)\cdot a_i=\max_{1\le i\le r} a_i
\cdot\sum_{i=1}^r a_i\le\max_{1\le i\le r} a_i\,.
}
$$
\eproof

Now we can prove the required bounds on the cardinalities $|B_n^d|$, $n\in\N$.

\begin{prop}\label{prop_posl}
Let $d\in\N$. Then the following holds.
\begin{enumerate}
\item Let $d\le2$. Then there exists a~constant $c>0$ such that 
$$
(2d)^{-n}|B_n^d|\ge cn^{-d/2}\,\text{ for every even }\,n\in\N\,.
$$
\item  Let $d\ge3$. Then there exists a~constant $c=c(d)>0$ such that 
$$
(2d)^{-n}|B_n^d|\le cn^{-d/2}\,\text{ for every }\,n\in\N\,.
$$
\end{enumerate}
\end{prop}
\proof
1. This follows from Proposition~\ref{prop_dMost2} and Lemma~\ref{lem_Stirling}.

2. We assume that $n=2m$ is even and compute the bound (which trivially holds for odd $n$):
\begin{eqnarray*}
(2d)^{-n}|B_n^d|&=&{\textstyle
(2d)^{-2m}\sum_{\substack{j_1,\,\ds,\,j_d\in\N_0\\j_1+\ds+j_d=m}}
\binom{2m}{j_1,\,j_1,\,
j_2,\,j_2,\,\ds,\,j_d,\,j_d}}\\
&=&{\textstyle
(2d)^{-2m}d^{2m}\binom{2m}{m}\sum_{\cdots}
\big(d^{-m}\binom{m}{j_1,\,j_2,\,\,\ds,\,j_d}\big)^2}\\
&\le&{\textstyle
4^{-m}\binom{2m}{m}d^{-m}\max_{\cdots}
\binom{m}{j_1,\,j_2,\,\,\ds,\,j_d}}\\
&=&{\textstyle
4^{-m}\binom{2m}{m}\cdot d^{-m}
\binom{m}{m_1,\,m_2,\,\,\ds,\,m_d}\,\ (\text{with }\, 
m_i\in\{\lfloor\frac{m}{d}\rfloor,\lceil\frac{m}{d}\rceil\})}\\
&\ll&m^{-1/2}\cdot m^{-(d-1)/2}=m^{-d/2}\ll n^{-d/2}\,\text{ for every even }\,n\in\N\,. 
\end{eqnarray*}
The first equality follows from Proposition~\ref{prop_B_n}. The 
second one is an algebraic rearrangement. The third inequality 
follows from Lemma~\ref{lem_treti} and from the multinomial expansion of 
$d^m=(1+1+\ds+1)^m$. The fourth equality follows from 
Lemma~\ref{lem_maxim}. The first factor in the fifth asymptotic 
inequality  follows from 
Lemma~\ref{lem_Stirling}. So does the second factor, but now we provide 
details; in general, $m_i\ne\frac{m}{d}$ and one cannot simply plug $\frac{m}{d}$ into Stirling's formula. 

Since $m_i=\frac{m}{d}+O(1)$ for every $m\in\N$, first
\begin{eqnarray*}
{\textstyle
\prod_{i=1}^d\big(\frac{m}{m_i}\big)^{m_i}}&=&
{\textstyle
\exp\big(
\sum_{i=1}^d m_i\log(\frac{m}{m_i})
\big)}\\
&=&{\textstyle
\exp\big(\sum_{i=1}^d
(\frac{m}{d}+O(1))\log(d+O(m^{-1}))\big)}\\
&=&{\textstyle
\exp\big(\sum_{i=1}^d
(\frac{m}{d}+O(1))(\log d+O(m^{-1}))\big)}\\
&=&\exp(m\log d+O(1))\ll d^m\,\text{ for every }\,m\in\N\,.
\end{eqnarray*}
Since $m_1+\ds+m_d=m$, we indeed compute with the help of Lemma~\ref{lem_Stirling} and the previous bound that ($c_1=\sqrt{2\pi}$ and $c_2>0$)
\begin{eqnarray*}
&&{\textstyle
d^{-m}
\binom{m}{m_1,\,m_2,\,\,\ds,\,m_d}   =d^{-m}
\frac{m!}{m_1!\cdot m_2!\cdot\ldots\cdot m_d!}}\\
&&{\textstyle
=d^{-m} 
c_1(1+O(m^{-1}))m^{1/2}\big(\frac{m}{\mathrm{e}}\big)^m/
\prod_{i=1}^d 
c_1(1+O(m_i^{-1}))m_i^{1/2}\big(\frac{m_i}{\mathrm{e}}\big)^{m_i}}\\
&&{\textstyle
=d^{-m}c_1^{1-d}(1+O(m^{-1}))
m^{1/2}\prod_{i=1}^d m_i^{-1/2}m^m
\prod_{i=1}^d m_i^{-m_i}}\\
&&{\textstyle
=d^{-m}c_2(1+O(m^{-1}))\cdot m^{(1-d)/2}
(1+O(m^{-1}))\cdot
\prod_{i=1}^d\big(\frac{m}{m_i}\big)^{m_i}}\\
&&\ll d^{-m}\cdot m^{(1-d)/2}\cdot d^m=m^{-(d-1)/2}\,\text{ for 
every }\,m\in\N\,.
\end{eqnarray*}

We can finish the justification of the original computation: the sixth and 
seventh steps are just algebraic rearrangements. 
\eproof

\noindent
We are done. The required bounds on
$${\textstyle
F_B^d(1)=\sum_{n\ge0}h(B_n^d)1^n=\sum_{n\ge0}(2d)^{-n}|B_n^d|
}
$$
follow from Proposition~\ref{prop_posl}, from the well known result that the series 
$${\textstyle
\sum_{n\ge1}n^{-\al}
}
$$ 
(absolutely) converges iff $\al>1$, and from $|B_0^d|=1$ and 
$|B_2^d|\ge1$.

\end{document}